\documentclass[11pt,a4paper]{article}

\usepackage{amsmath}     
\usepackage{dsfont}		
\usepackage{amsfonts}		
\usepackage{amsthm}
\usepackage{graphicx}
\usepackage{amssymb}

\usepackage{latexsym}
\usepackage{euscript,makeidx,color,mathrsfs}
\usepackage{enumerate}
\usepackage{tabu}

\usepackage[colorlinks,linkcolor=blue,anchorcolor=green,citecolor=red]{hyperref}

\def\sqr#1#2{{\vcenter{\vbox{\hrule height.#2pt
              \hbox{\vrule width.#2pt height#1pt \kern#1pt \vrule width.#2pt}
              \hrule height.#2pt}}}}
\def\5n{\negthinspace \negthinspace \negthinspace \negthinspace \negthinspace }
\def\4n{\negthinspace \negthinspace \negthinspace \negthinspace }
\def\3n{\negthinspace \negthinspace \negthinspace }
\def\2n{\negthinspace \negthinspace }
\def\1n{\negthinspace }

\def\dbD{\mathbb{D}}
\def\dbE{\mathbb{E}}
\def\dbF{\mathbb{F}}

\def\dbH{\mathbb{H}}

\def\dbK{\mathbb{K}}
\def\dbL{\mathbb{L}}

\def\dbN{\mathbb{N}}

\def\dbR{\mathbb{R}}

\def\dbU{\mathbb{U}}
\def\dbV{\mathbb{V}}
\def\dbW{\mathbb{W}}
\def\dbX{\mathbb{X}}

\def\={\buildrel \triangle \over =}

\def\ds{\displaystyle}

\def\nb{\noalign{\bs}}

\def\a{\alpha}

\def\d{\delta}
\def\e{\varepsilon}

\def\si{\sigma}
\def\t{\tau}
\def\f{\varphi}
\def\th{\theta}

\def\G{\Gamma}
\def\D{\Delta}

\def\O{\Omega}
\def\mf{\mathcal{F}}
\def\me{\mathbb{E}}
\def\rd{\,\mathrm d}  
\def\bal{\begin{aligned}}
\def\eal{\end{aligned}}
\def\nb{\nabla}
\def\cF{{\cal F}}

\def\cJ{{\cal J}}

\def\cL{{\cal L}}

\def\cP{{\cal P}}

\def\cR{{\cal R}}
\def\cS{{\cal S}}
\def\cT{{\cal T}}

\def\no{\noindent}

\def\ms{\medskip}

\def\q{\quad}
\def\qq{\qquad}

\def\lt{\left}
\def\rt{\right}

\def\rf{\eqref}

\def\h{\widehat}
\def\wt{\widetilde}

\def\cd{\cdot}
\def\cds{\cdots}

\def\ae{\hbox{\rm a.e.}}

\def\({\Big (}
\def\){\Big )}
\def\[{\Big[}
\def\]{\Big]}

\def\bde{\begin{definition}\label}
\def\ede{\end{definition}}
\def\be{\begin{equation}}
\def\bel{\begin{equation}\label}
\def\ee{\end{equation}}
\def\beq{\begin{equation*}\begin{aligned}}
\def\eeq{\end{aligned}\end{equation*}}
\def\bt{\begin{theorem}\label}
\def\et{\end{theorem}}
\def\bc{\begin{corollary}\label}
\def\ec{\end{corollary}}
\def\bl{\begin{lemma}\label}
\def\el{\end{lemma}}
\def\bp{\begin{proposition}\label}
\def\ep{\end{proposition}}
\def\bas{\begin{assumption}\label}
\def\eas{\end{assumption}}
\def\br{\begin{remark}\label}
\def\er{\end{remark}}
\def\bex{\begin{example}\label}
\def\ex{\end{example}}
\def\ba{\begin{array}}
\def\ea{\end{array}}
\def\ed{\end{document}}

\def\square#1{\vbox{\hrule\hbox{\vrule height#1%
     \kern#1\vrule}\hrule}}
\def\rectangle#1#2{\vbox{\hrule\hbox{\vrule height#1%
     \kern#2\vrule}\hrule}}

\font\tenbb=msbm10 \font\sevenbb=msbm7 \font\fivebb=msbm5

\newfam\bbfam
\scriptscriptfont\bbfam=\fivebb \textfont\bbfam=\tenbb
\scriptfont\bbfam=\sevenbb

\newtheorem{theorem}{\hskip 1.3em Theorem}[section]
\newtheorem{definition}[theorem]{\hskip 1.3em Definition}
\newtheorem{proposition}[theorem]{\hskip 1.3em Proposition}
\newtheorem{corollary}[theorem]{\hskip 1.3em Corollary}
\newtheorem{lemma}[theorem]{\hskip 1.3em Lemma}
\newtheorem{remark}[theorem]{\hskip 1.3em Remark}
\newtheorem{example}[theorem]{\hskip 1.3em Example}
\newtheorem{algorithm}[theorem]{\hskip 1.3em Algorithm}

\newtheorem{assumption}[theorem]{\hskip 1.3em Assumption}

\makeatletter
   
   \@addtoreset{equation}{section}
\makeatother

\begin{document}

\title{Strong Error Estimates for a Space-Time Discretization of the
Linear-Quadratic Control Problem with the Stochastic Heat Equation with Linear Noise\thanks{This work is supported in part by
the National Natural Science Foundation of China (11801467), and the Chongqing Natural Science Foundation
(cstc2018jcyjAX0148).}}

\author{Andreas Prohl\thanks{
Mathematisches Institut, Universit\"at T\"ubingen, Auf der Morgenstelle 10,
D-72076 T\"ubingen, Germany.
 {\small\it
e-mail:} {\small\tt prohl@na.uni-tuebingen.de}.} \quad 
and \quad Yanqing Wang\thanks{Corresponding author.
School of Mathematics and Statistics, Southwest University, Chongqing 400715, China.  {\small\it
e-mail:} {\small\tt yqwang@amss.ac.cn}. \ms}}

\date{November 28, 2020}
\maketitle

\begin{abstract}
We propose a time-implicit, finite-element based space-time discretization
of  the necessary and sufficient optimality conditions for
the stochastic linear-quadratic optimal control problem 
with the stochastic heat equation driven by {\em linear} noise of type $[X(t)+\sigma(t)]{\rm d}W(t)$, and prove
{optimal} convergence {\em w.r.t.} both, space and time discretization parameters.
In particular, we employ the stochastic Riccati equation as a proper analytical tool to handle the linear noise, and thus extend the applicability of the earlier work \cite{Prohl-Wang20}, where the error analysis was restricted to 
{\em additive} noise.
%
%
%
%
%
%
%
%
\end{abstract}

\ms

\no\bf Keywords: \rm Error estimate, stochastic linear quadratic problem, 
 stochastic heat equation,
Pontryagin's maximum principle,
stochastic Riccati equation


\ms

\no\bf AMS 2010 subject classification: \rm 49J20,
 65M60,
 93E20

\section{Introduction}

Let $D \subset {\mathbb R}^d$ be a bounded domain with $C^2$ boundary and $T>0$ be given. 
Our goal is to numerically approximate the ${\mathbb F}$-adapted control process 
$U^* \equiv \{ U^*(t);\, t \in [0,T]\}$ on the filtered probability space 
$(\Omega, {\mathcal F}, {\mathbb F}, {\mathbb P})$ that minimizes the quadratic  functional ($\alpha \geq 0$)
\begin{equation} \label{w1003e2}
{\mathcal J}(X, U)=\frac 1 2 {\mathbb E} \Bigl[\int_0^T  \Vert X(t)  \Vert_{\dbL^2}^2+ \Vert U(t) \Vert^2_{\dbL^2} \, {\rm d}t
+ \alpha \Vert X(T)  \Vert^2_{\dbL^2}\Bigr] 
\end{equation}
subject to the (controlled forward) stochastic PDE ({\bf SPDE}, for short) of the form
\bel{w1013e1}
\lt\{
\begin{array}{ll}
{\rm d}X(t)=\bigl[\D X(t)+U(t)\bigr]\, {\rm d}t+\big[X(t)+\si(t)\big] {\rm d}W(t) \q &\forall\, t \in [0,T]\,,\\
X(0)=X_0\,,
\end{array}
\rt.
\ee
for a proper function $\sigma: [0,T] \rightarrow {\mathbb L}^2$, a suitable initial datum $X_0$,
homogeneous Dirichlet boundary conditions, and a Wiener process
$W:= \{ W(t):\, t \in [0,T]\}$, which here is ${\mathbb R}$-valued for the sake of simplicity.
A unique (strong) minimizer $(X^*, U^*) $ 
of the stochastic optimal control problem: `minimize \eqref{w1003e2} subject to \eqref{w1013e1}' may then be deduced --- which we below refer to as problem {\bf SLQ}; see e.g.~\cite{Bensoussan83,Lv-Zhang20}, and {Section \ref{pre} for a further specification of the data.}

{\bf SLQ} is a prototypic stochastic optimization problem on (infinite-dimensional) Hilbert spaces, for which the numerical analysis so-far is rare in the literature; {\em cf.}~the references in \cite{Prohl-Wang20}. In the foregoing work \cite{Prohl-Wang20}, {optimal} strong error estimates were shown for a space-time discretization of a corresponding problem where the equation \rf{w1013e1}$_1$ was driven by {\em additive} noise. On the level of Pontryagin's maximum principle that we apply below to develop numerical methods for {\bf SLQ} this difference, in particular, simplified a lot the (numerical analysis of the) backward 
stochastic PDE ({\bf BSPDE}, for short), which for the current problem {\bf SLQ} 
with solution tuple $(Y,Z)$ reads:
\bel{bshe1}
\lt\{
\begin{array}{ll}
{\rm d}Y(t)= \bigl[-\D Y(t)-Z(t)+ X(t)  \bigr]{\rm d}t+Z(t) {\rm d}W(t)  \q & \forall\, t \in [0,T]\, ,\\
Y(T)=- \a X(T)\, .
\end{array}
\rt.
\ee
Note that in the case of {\em additive} noise in 
\rf{w1013e1},  $Z(t)$ in
\rf{bshe1} does not appear in the drift term, which is the reason why the tools
which were developed for the corresponding numerical analysis in \cite{Prohl-Wang20} do to cover 
the case of {\em linear} noise as present in \rf{w1013e1}.
 For {\bf SLQ}, the optimality system consists of  \rf{w1013e1}, \rf{bshe1} and 
 Pontryagin's maximum {condition}\begin{equation}\label{pontr1} 
0 = U(t) - Y(t) \qquad \forall\, t \in (0,T)\, ,
\end{equation}
which then uniquely determines the optimal process tuple denoted by $(X^*,U^*)$ of problem {\bf SLQ}.

Based on Pontryagin's maximum principle, problem {\bf SLQ}   may be numerically accessed by solving the coupled system \rf{w1013e1}, \rf{bshe1} and \rf{pontr1}, which is a 
forward-backward stochastic PDEs ({\bf FBSPDE}, for short). 
In order to derive strong error estimates in \cite{Prohl-Wang20} for the related stochastic control problem with {\em additive} noise in \rf{w1013e1},  a spatial semi-discretization via finite elements (with step size parameter $h$ for the mesh) with solution $(X^*_h, Y_h, Z_h, U^*_h)$ was considered in a first step, for which optimal convergence rates were obtained; the key link to show optimal strong error estimates for 
the full space-time discretization (with additional time step parameter $\tau$) in a second step then depended on the verification of $h$-independent stability results for the above solution quadruple $(X^*_h, Y_h, Z_h, U^*_h)$, whose derivation via Malliavin calculus rested on the fact that $Z_h$ did not enter the drift part in  (the finite element version of) the corresponding modification of {\bf BSPDE} \rf{bshe1} within (the finite element version of) the coupled optimality system {\bf FBSPDE}.
In this work, we use the stochastic Riccati equation \rf{riccati} as proper representation tool for the solution of the semi-discretization {\bf SLQ}$_h$ (see \rf{w1003e2h}--\rf{w1013e1a})   to deduce the relevant
$h$-independent stability results for the solution quadruple $(X^*_h, Y_h, Z_h, U^*_h)$ of the related optimality system ({\bf FBSPDE})$_h$, which is \rf{w1013e1a}--\rf{bshe1a}; these results may then be used to prove optimal convergence rates for optimal tuple of the space-time discretization {\bf SLQ}$_{h\t}$, which is
\rf{w1212e3}--\rf{w1003e12} in equivalent form. Specifically, a relevant result
 is the following, which bounds the temporal variation of the component $Z_h$ of the solution to ({\bf FBSPDE})$_h$ (see Lemmata \ref{w229l3}, \ref{w1024l1}),
\begin{equation}\label{z19}\me\bigl[\|Z_h(t)-Z_h(s)\|_{\dbL^2}^2\bigr]\leq C|t-s|\qq \forall \, t,\,s \in [0,T]\,,
\end{equation}
where {\em $C>0$ is independent of   $h$}.
We remark that \rf{w1212e3}--\rf{w1003e12} is a modification of the implicit Euler method, which is again due to the role that $Z_h$ plays in the drift part in  ({\bf FBSPDE})$_h$, and  that needs be properly addressed numerically; see Remark \ref{w829r1}.
{\color{cyan}
}

\ms

To computationally solve this discrete, coupled optimality system requires huge computational resources; instead,
we again return to the fully discretized problem {\bf SLQ}$_{h\t}$
\eqref{w1003e8}--\eqref{w1003e7}
 and exploit its character as a minimization problem to initiate a {\em decoupled} gradient descent method to
successively determine approximations of the optimal control; this method, which is close to the one in \cite{Dunst-Prohl16} where also computational experiments are provided, is detailed in
Section \ref{numopt}, and {an optimal} convergence rate is shown for this iteration --- which is the second goal in this work. 

\ms

The rest of this paper is organized as follows. In Section \ref{pre}, we introduce  notations, and review relevant rates
of convergence for a discretization in space and time of {\bf SPDE} \eqref{w1013e1} and a semi-discretization in space
of  {\bf BSPDE} \eqref{bshe1}, which are both needed in 
 Section \ref{rate}. In Section \ref{rate}, we prove  rates of convergence for a space-time discretization of a coupled {\bf FBSPDE}, which is related to problem {\bf SLQ}. Convergence of the related iterative gradient descent method towards the
optimal pair $(X^*,U^*)$ of  problem {\bf SLQ} is shown in Section \ref{numopt}.

\section{Preliminaries}\label{pre}

\subsection{Notations {and assumptions} --- involved processes and the finite element method}\label{not1}

Let $\bigl({\mathbb K}, ( \cdot, \cdot)_{{\mathbb K}}\bigr)$ be a separable Hilbert space. 
By $\Vert \cdot \Vert_{{\mathbb L}^2}$ resp.~$(\cdot, \cdot)_{\dbL^2}$, 
we denote the norm resp.~the scalar product in Lebesgue space ${\mathbb L}^2 := L^2(D)$.
By $\|\cd\|_{\dbH_0^1}$, $\|\cd\|_{\dbH^2}$, $\|\cd\|_{\dbH^3}$, we denote norms in Sobolev spaces ${\mathbb H}^1_0:=H_0^1(D)$, $\dbH^2:=H^2(D)$, $\dbH^3:=H^3(D)$ respectively.
Let $(\Omega, {\mathcal F}, {\mathbb F}, {\mathbb P})$ be a complete filtered probability space, where 
${\mathbb F}=\{\mf_t\}_{t\in[0,T]}$ is the filtration generated by the ${\mathbb R}$-valued Wiener process $W$, which is augmented by 
all the ${\mathbb P}$-null sets. The space of all ${\mathbb F}$-adapted processes 
$X: \Omega \times [0,T] \rightarrow {\mathbb K}$ satisfying 
${\mathbb E}[\int_0^T \Vert X(t)\Vert^2_{{\mathbb K}}\, {\rm d}t] < \infty$ is denoted by 
{$L^2_\dbF(0,T; {\mathbb K})$}; the space of all ${\mathbb F}$-adapted processes 
$X: \Omega \times [0,T] \rightarrow {\mathbb K}$ with continuous path satisfying 
${\mathbb E}[\sup_{t \in [0,T]} \Vert X(t)\Vert^2_{\mathbb K}] < \infty$ is denoted by
$L^2_{{\mathbb F}}\bigl(\Omega; C([0,T]; {\mathbb K})\bigr)$;
for any $t\in [0,T]$, the space of $\dbK$-valued $\mf_t$-measurable random variables $\eta$ satisfying $\me[\|\eta\|^2_{\dbK}]<\infty$
is denoted by $L^2_{\mf_t}(\O;\dbK)$.

\ms

We partition the bounded domain $D \subset {\mathbb R}^d$ via a regular triangulation ${\mathcal T}_h$ into elements $K$ with maximum mesh size
$h := \max \{ {\rm diam}(K):\, K \in {\mathcal T}_h\}$, and consider space
\beq
{\mathbb V}_h :=  \lt\{\phi \in {\mathbb H}^1_0:\, \phi \bigl\vert_K \in {\mathbb P}_1(K) \quad \forall\, K \in {\mathcal T}_h \rt\}\,, 
\eeq
where ${\mathbb P}_1(K)$ denotes the space of {polynomials} of degree $1$; see {\em e.g.}~\cite{Brenner-Scott08}.
We define the discrete Laplacian $\Delta_h: {\mathbb V}_h \rightarrow {\mathbb V}_h$ by $(-\Delta_h \xi_h, \phi_h)_{{\mathbb L}^2} = (\nabla \xi_h, \nabla \phi_h)_{{\mathbb L}^2}$ for all $\xi_h, \phi_h \in {\mathbb V}_h$, the ${\mathbb L}^2$-projection $\Pi_h: {\mathbb L}^2 \rightarrow {\mathbb V}_h$ by $(\Pi_h \xi - \xi, \phi_h)_{{\mathbb L}^2}= 0$ for all $\xi\in \dbL^2,\, \phi_h \in {\mathbb V}_h$, and the Ritz-projection $\cR_h:\dbH_0^1\to \dbV_h$
by $(\nb[\cR_h\xi-\xi], \nb\phi_h)_{\dbL^2}=0$ for all $\xi\in\dbH_0^1,\,\phi_h\in \dbV_h$.  
Via definition of $\cR_h$, it is easy to get that
\bel{w1024e1}
\|\nb \cR_h\xi\|_{\dbL^2}\leq \|\nb\xi\|_{\dbL^2}\qq \forall \, \xi \in \dbH_0^1\,.
\ee


We denote by $I_\tau = \{ t_{n}\}_{n=0}^N \subset [0,T]$ a time mesh with maximum step size $\tau :=\max\{t_{n+1}-t_n:\,n=0,1,\cds,N-1\}$, and
$\D_nW:=W(t_n)-W(t_{n-1})$ for all $n=1,\cds,N$. For given time mesh $I_\t$, we can define $\mu(\cd),\,\nu(\cd)$ by
\bel{w827e1}
\mu(t):=t_{n+1}  \, ,\qquad  \nu(t):=t_n  \qquad \forall\, t \in [t_n, t_{n+1}),\,n=0,1,\cds,N-1\, .
\ee
Throughout this work, we assume that $\t\leq 1$.
For simplicity, we choose a uniform partition, i.e. $\t=T/N$. 
The results in this work still hold for quasi-uniform partitions. 
Throughout this work, we shall make use of the following\\
{\bf Assumption (A)}: $X_0\in  \dbH_0^1\cap\dbH^3$, and $\si\in C^1(0,T;\dbH_0^1\cap\dbH^3).$

\subsection{The stochastic heat equation and rates of strong convergence for its space-time discretization}\label{opt-1}

Under Assumption {\bf(A)}, in particular,
and given 
$U\in L^2_{\mathbb F} (0,T; \dbH_0^1 )$
in
{\bf SPDE} \eqref{w1013e1}, there exists a unique strong solution $X \in L^2_{{\mathbb F}}\bigl(\Omega; C([0,T]; {\mathbb H}^1_0)\bigr) \cap L^2_{{\mathbb F}}\bigl(0,T;  {\mathbb H}^2\bigr)$
which satisfies  
the following estimate (see {\em e.g.}~\cite{Chow15}),
\begin{equation}\label{stochheat1} 
{\mathbb E}\Big[ \sup_{t \in [0,T]} \Vert X(t)\Vert_{{\mathbb H}_0^1}^2 +
\int_0^T \Vert X(t)\Vert^2_{{\mathbb H}^2}\, {\rm d}t\Big] 
\leq C  {\mathbb E}\Big[ \Vert X_0\Vert^2_{{\mathbb H}_0^1} + \int_0^T \| U(t)\|^2_{\dbL^2}+\| \si(t)\|^2_{\dbH_0^1}\, {\rm d}t
\Big]\, .
\end{equation}
Obviously, $X$ satisfies the following variational form ${\mathbb P}$-a.s.~for all $t \in [0,T]$,
\begin{equation}\label{forw1}
\bal
&\big( X(t), \phi\big)_{\dbL^2} - \big(X_0, \phi\big)_{\dbL^2} + \int_0^t  \big(\nabla X(s), \nabla \phi\big)_{\dbL^2} 
- \big(U(s), \phi\big)_{\dbL^2} \, {\rm d}s\\
&\qq\qq= \int_0^t \big( X(s)+\si(s), \phi \big)_{\dbL^2} \rd W(s)
\qquad \forall\, \phi \in {\mathbb H}^1_0\, .
\eal
\end{equation}

A finite element discretization of \eqref{forw1}  which we later refer to as {\bf SPDE}$_h$  then reads: For all $t \in [0,T]$, find
$X_h \in L^2_{{\mathbb F}}\bigl(\Omega; C([0,T]; {\mathbb V}_h)\bigr)$ such that ${\mathbb P}$-a.s.~and for all times $t \in [0,T]$
\bel{spdedisch1}
\bal
&\big(X_h(t), \phi_h\big)_{\dbL^2} - \big(X_h(0), \phi_h\big)_{\dbL^2} + \int_0^t  \big(\nabla X_h(s), \nabla \phi_h\big)_{\dbL^2} 
- \big( U(s), \phi_h\big)_{\dbL^2} \rd s \\ 
&\qquad \qq = \int_0^t \big( X_h(s)+\cR_h\si(s), \phi_h \big)_{\dbL^2} \rd W(s), \qquad \forall\, \phi_h \in {\mathbb V}_h \, .
\eal
\ee
Equation \eqref{spdedisch1} may be recast into the following stochastic differential equation,
\bel{sde}
\lt\{
\bal
&dX_h(t)=\bigl[\Delta_hX_h(t)+\Pi_hU(t)\bigr] {\rm d}t+ \big[X_h(t)+\cR_h\si(t)\big] {\rm d}W(t)\q \forall\, t\in [0,T]\, , \\
&X_h(0)= \cR_hX_{0}\, .
\eal
\rt.
\ee
{Thanks to this equivalence,  we do not distinguish between {\bf SPDE}$_h$ \eqref{spdedisch1} and equation \eqref{sde} 
throughout this paper.}

The derivation of an error estimate is well-known {(see {\em e.g.}~\cite{Yan05})}, which uses the improved (spatial) regularity properties of the strong solution to deduce
\begin{equation}\label{esti-space1}
\sup_{t \in [0,T]} {\mathbb E} \bigl[\Vert X_h(t) - X(t)\Vert^2_{{\mathbb L}^2}\bigr]
+ {\mathbb E} \Bigl[\int_0^T  \Vert \nabla \bigl[ X_h(t) - X(t)\bigr]\Vert^2_{{\mathbb L}^2}\, {\rm d}t\Bigr]
\leq C h^2\, .
\end{equation}

We now consider a time-implicit discretization of \eqref{spdedisch1} on a partition $I_\tau$ of $[0,T]$. The problem then reads: For every $0 \leq n \leq N-1$, find a solution $X^{n+1}_h \in L^2_{{\mathcal F}_{t_{n+1}}}(\Omega; {\mathbb V}_h)$ such that ${\mathbb P}$-a.s.
\begin{equation}\label{esti-time1} 
\big(X^{n+1}_h - X^{n}_h, \phi_h\big)_{\dbL^2} + \tau \Big[\big(\nabla X_h^{n+1},\nabla \phi_h\big)_{\dbL^2} 
- \bigl(U(t_n), \phi_h\bigr)_{\dbL^2}\Big] 
= \big( X_h^n, \phi_h\big)_{\dbL^2}\Delta_{n+1} W\, .
\end{equation}
The verification of the error estimate (see \cite{Yan05})
\begin{equation}\label{euler1}
\max_{0 \leq n \leq N} {\mathbb E} \bigl[\Vert X_h(t_n) - X^n_h\Vert^2_{{\mathbb L}^2} \bigr]
+ \tau \sum_{n=1}^N  {\mathbb E} \Bigl[ \Vert \nabla \bigl[ X_h(t_n) - X^n_h\bigr] \bigr\Vert^2_{{\mathbb L}^2}\Bigr]
\leq C \tau
\end{equation}
rests on 
 stability properties of the implicit Euler, Assumption {\bf (A)}, as well as the further assumption that
\beq
\sum_{n=0}^{N-1} \me \Bigl[  \int_{t_n}^{t_{n+1}}\|U(t)-U(t_n)\|_{\dbL^2}^2 \rd t \Bigr]  \leq C \t\, .
\eeq

{When solving problem {\bf SLQ}, we select the
optimal control $U^*$ for $U$, in particular, which inherits improved regularity conditions from the solution of the adjoint equation to the optimal control $U^*$ via the optimality condition.
}


\subsection{The backward stochastic heat equation --- a finite element based spatial discretization}\label{opt-2}
Let $Y_T \in L^2_{{\mathcal F}_T}\bigl( \Omega; {\mathbb H}^1_0\bigr)$ and $f \in L^2_{{\mathbb F}} (0,T; {\mathbb L}^2)$. A strong solution to the backward stochastic heat equation
\bel{bshe}
\lt\{
\begin{array}{ll}
{\rm d}Y(t)= \bigl[-\D Y(t)-Z(t)+f(t) \bigr] \rd t+Z(t) \rd W(t)  & \forall\, t \in [0,T]\, , \\
Y(T)=Y_T  
\end{array}
\rt.
\ee
{with homogeneous Dirichlet boundary data}
is a pair of square integrable ${\mathbb F}$-adapted processes 
$(Y,Z) \in \bigl(L^2_{{\mathbb F}}\bigl(\Omega; C([0,T]; {\mathbb H}^1_0) \cap L^2_\dbF\bigl(0,T; {\mathbb H}^1_0 \cap {\mathbb H}^2\bigr) \bigr)\times L^2_\dbF\bigl(0,T; {\mathbb H}^1_0\bigr)$ that satisfies the following variational form
 ${\mathbb P}$-a.s.~for all  $t \in [0,T]$,
 \begin{equation}\label{vari-1}
 \bal
&\bigl( Y_T, \phi\bigr)_{\dbL^2} - \big(Y(t),\phi\big)_{\dbL^2} - \int_t^T  \bigl(\nabla Y(s), \nabla \phi \bigr)_{\dbL^2} 
- \big( Z(s), \phi\big)_{\dbL^2}
+  \bigl( f(s), \phi\bigr)_{\dbL^2} \, {\rm d}s \\
&\qq\qq= \int_t^T \bigl( Z(s), \phi\bigr)_{\dbL^2}\, {\rm d}W(s)
 \qquad \forall\, \phi \in {\mathbb H}^1_0\,.
 \eal
 \end{equation}
 The existence of a strong  solution to \eqref{bshe}, as well as its uniqueness are shown in \cite{Du-Tang12}; moreover, 
there exists a constant $C \equiv C(D,T) > 0$ such that
 \begin{equation}\label{vari-2}
 {\mathbb E}\bigl[ \sup_{t \in [0,T]} \Vert Y(t)\Vert^2_{\dbH_0^1}\bigr]
 + {\mathbb E}\Bigl[\int_0^T  \Vert Y(t)\Vert^2_{{\mathbb H}^2} + 
 \Vert Z(t) \Vert^2_{\dbH_0^1}  \rd t\Bigr]
  \leq
 C {\mathbb E}\Bigl[  \Vert Y_T\Vert^2_{\dbH_0^1} +  \int_0^T \Vert f(t) \Vert^2_{{\mathbb L}^2} \rd t\Bigr]\, .
 \end{equation}
 
We may consider a finite element discretization  of the {\bf BSPDE} \eqref{bshe}.
Let $Y_{T,h}=  \cR_h Y_T \in L^2_{{\mathcal F}_T}(\Omega; {\mathbb V}_h)$, 
 an approximation of $Y_T$. The problem {\bf BSPDE}$_h$ then reads:
 Find a pair $(Y_h, Z_h) \in L^2_{{\mathbb F}}\bigl( \Omega; C([0,T]; {\mathbb V}_h)\bigr) \times L^2_{{\mathbb F}}\bigl( 0,T; {\mathbb V}_h\bigr)$ such that ${\mathbb P}$-a.s.~for all $t \in [0,T]$ 
\bel{vf2}
\bal
\big(Y_{T,h},\phi_h\big)_{\dbL^2} - \big( Y_h(t), \phi_h\big)_{\dbL^2}  
=& \int_t^T  \big(\nabla Y_h(s), \nabla \phi_h \big)_{\dbL^2}
- \big( Z(s), \phi_h\big)_{\dbL^2}+  \big( f(s), \phi_h\big)_{\dbL^2} \, {\rm d}s \\ 
 &+ \int_t^T \big( Z_h(s), \phi_h\big)_{\dbL^2}\, {\rm d}W(s)
 \qquad \forall\, \phi_h \in {\mathbb V}_h\,.
\eal
\ee
Actually, equation \eqref{vf2} is equivalent to the following BSDE:
\bel{bsde}
\lt\{
\bal
&{\rm d}Y_h(t)=\bigl[-\Delta_hY_h(t)-Z_h(t)+ \Pi_h f(t)\bigr] {\rm d}t+ Z_h(t){\rm d}W(t) \quad \forall\, t\in [0,T]\, ,\\
&Y_h(T)=Y_{T,h}\, ,
\eal
\rt.
\ee
and we do not distinguish between {\bf BSPDE}$_h$ \eqref{vf2} and BSDE \eqref{bsde} 
throughout this paper.
The existence and uniqueness of a solution tuple $(Y_h, Z_h)$ follows from \cite[Theorem 2.1]{ElKaroui-Peng-Quenez97}. Moreover, there exists $C \equiv C( T) >0$ such that
\beq
\bal
\sup_{t \in [0,T]} {\mathbb E} \bigl[\Vert \nabla Y_h(t)\Vert^2_{{\mathbb L}^2}\bigr] +
{\mathbb E} \Bigl[\int_0^T  \Vert \Delta_h Y_h(t)\Vert^2_{{\mathbb L}^2} +
\Vert \nabla Z_h(t)\Vert^2_{{\mathbb L}^2} \, {\rm d}t\Bigr]
\leq
C {\mathbb E}\Bigl[ \lt\| \nabla Y_{T,h}\rt\|^2_{{\mathbb L}^2}+ \int_0^T\|f(t)\|^2_{\dbL^2}\rd t\Bigr]\, ;
\eal
\eeq
cf.~\cite[Lemma 3.1]{Dunst-Prohl16}. The  following result is taken from \cite[Theorem 3.2]{Dunst-Prohl16}, whose proof exploits  the bounds \eqref{vari-2}.

\begin{theorem}\label{w0911t1}
Assume $Y_T \in L^2_{{\mathcal F}_T}(\Omega; {\mathbb H}^1_0)$ and $f\in L^2_\dbF(0,T;\dbL^2)$. Let $(Y,Z)$ be the solution of \eqref{vari-1}, and $(Y_h, Z_h)$ solves \eqref{vf2}. There exists
$C \equiv C(Y_T, f,T) >0$ such that
\beq
\sup_{t \in [0,T]} \me \bigl[\|Y(t)-Y_h(t)\|^2_{{\mathbb L}^2} \bigr]
+\me \Bigl[\int_0^T  \|\nabla\bigl[ Y(t)-Y_h(t)\bigr]\|^2_{{\mathbb L}^2}+\|Z(t)-Z_h(t)\|^2_{{\mathbb L}^2}\, {\rm d}t\Bigr]
 \leq C h^2\, .
\eeq
\end{theorem}

{When solving problem {\bf SLQ}, we may take $f=X^*$ in {\bf BSPDE} \rf{bshe}, 
and $\Pi_h f=X^*_h$ in \rf{bsde}, 
where $X^*, X^*_h$ are the optimal states of problems {\bf SLQ}, and {\bf SLQ}$_h$, respectively.
}


\ms
\subsection{Temporal discretization of problem {\bf SLQ} --- the role of Malliavin derivatives}\label{temp-mall}

The numerical analysis of a temporal discretization of problem {\bf SLQ} requires Malliavin calculus to bound temporal increments such as $\me[\Vert Z_h(t) - Z_h(s)\Vert^2_{{\mathbb L}^2}]$ in terms of $\vert t-s \vert$, where $s,t \in [0,T]$, 
and $Z_h$ is the second component of  BSDE \eqref{bsde} with $\Pi_h f=X_h^*$. We therefore briefly recall the definition and results of the Malliavin derivative
of processes, which will be applied below. For further details, we refer to \cite{Nualart06,ElKaroui-Peng-Quenez97}. 

We define the It\^o isometry $\dbW:L^2(0,T;{\mathbb R})\to L^2_{\cF_{T}}(\Omega;{\mathbb R})$ by
\beq
\dbW(g)=\int_0^{T}g(t) \, {\mathrm d}W(t)\, .
\eeq
For $\ell \in {\mathbb N}$, we denote by $C_p^\infty(\dbR^\ell)$ the {space} of all smooth functions $s:\dbR^\ell \to \dbR$ such that $s$ and all of its partial derivatives have  polynomial growth.
Let  $\cP$ be the {set} of ${\mathbb R}$-valued  random variables of the form
\beq
F=s\bigl(\dbW(g_1),\dbW(g_2),\cdots,\dbW(g_\ell) \bigr) 
\eeq
for some $s \in C_p^{\infty}(\dbR^\ell)$, $\ell \in {\mathbb N}$, and   $g_1,\ldots,g_\ell\in L^2(0,T;{\mathbb R})$.
To any $F\in \cP$ we define its  ${\mathbb R}$-valued Malliavin derivative 
$DF := \{ D_\th F;\, 0 \leq \th \leq T\}$ process via
\beq
D_\th F=\sum\limits_{i=1}^\ell \frac{\partial s}{\partial x_i} (\dbW(g_1),\dbW(g_2),\cdots,W(g_\ell))g_i(\th)\, .
\eeq
In general, we can define the $k$-th iterated derivative of $F$ by $D^kF=D(D^{k-1}F)$, for any
$k\in {\mathbb N}$. Note that, generally, for any $\theta \in[0,T]$, $D_\th F$ is $\mf_T$-measurable; and if $F$ is 
$\mf_t$-measurable, then $D_\th F=0$ for any $\th\in (t,T]$. 

Now we can extend the derivative operator to $\dbK$-valued variables. For any $k\in \dbN$, and $u$ in the set of $\dbK$-valued variables:
\beq
\cP_\dbK= \Bigl\{ u=\sum_{j=1}^n F_j \phi_j: F_j\in \cP,\,\phi_j\in \dbK,\, n\in \dbN  \Bigr\}\, ,
\eeq
we can define the $k$-th iterated derivative of $u$ by
$
D^k u=\sum_{j=1}^n D^kF_j\otimes \phi_j\, .
$
For $p \geq 1$, we define the norm $\|\cdot\|_{k,p}$  via
\beq
\|u\|_{k,p} :=\Big(\dbE \bigl[\|u\|_\dbK ^{p}\bigr] +\sum_{j=1}^k \me \bigl[\lt\|D^ju\rt\|_{\lt({ L^2(0,T;{\mathbb R})}\rt)^{\otimes j}\otimes \dbK}^p\bigr] \Big)^{1/ p}\, .
\eeq
Then $\dbD^{k,p}({\mathbb K})$ is the completion of $\cP_\dbK$ under the norm $\|\cdot\|_{k,p}$.

\section{Rates of convergence for a {spatio-temporal} discretization of problem SLQ}\label{rate}

\subsection{The discrete maximum principle and main results}\label{rate-1}

In this part, we discretize the original problem {\bf SLQ} within two steps, 
starting with its semi-discretization in space (which is referred to as {\bf SLQ}$_{h}$), which is then followed by a
discretization in space and time (which is referred to as  {\bf SLQ}$_{h\t}$). 
Our goal in this section is to prove rates of convergence in both
cases. By \cite{Lv-Zhang20}, problem {\bf SLQ} is uniquely solvable, and its optimal pair $(X^*,U^*)$ may be 
characterized by the following coupled {\bf FBSPDE}  {(supplemented by homogeneous Dirichlet data for $X^*, Y$}) with a unique solution
$(X^*, Y, Z, U^*)$, 
\bel{w212e3}
\lt\{
\bal
&{\rm d} X^*(t)=\big[\D X^*(t)+ U^*(t) \big]{\rm d}t
+  \big[X^*(t)+\si(t)\big] {\rm d}W(t)\qquad \forall\, t\in (0,T) \, ,\\
&{\rm d} Y(t)=\big[-\D Y(t)- Z(t)+ X^*(t) \big] {\rm d}t+Z(t){\rm d}W(t) \qquad \forall\,  t\in (0,T) \, ,\\
& X^*(0)=X_0\, ,\qquad Y(T)=-\a X^*(T)\, ,
\eal
\rt.
\ee
with the condition
\bel{w212e4}
 U^*(t)-Y(t)=0 \qq \forall\, t\in (0,T)\, .
\ee

By \eqref{w212e4}  and the fact that $Y\in L^2_\dbF(\O;C([0,T];\dbH_0^1))$, we find that the optimal control $U^*$ has  continuous paths, {taking zero values on the boundary $\partial D$}.

\ms
The spatial semi-discretization  {\bf SLQ}$_h$ of  problem {\bf SLQ} reads as follows:  Find an optimal pair
$(X_h^*, U^*_h) \in L^2_{{\mathbb F}}(\Omega; C([0,T]; {\mathbb V}_h)) \times L^2_\dbF(0,T; {\mathbb V}_h) $ that minimizes the functional
\begin{equation} \label{w1003e2h}
\cJ (X_h, U_h)=\frac 1 2 \me \Bigl[\int_0^T  \| X_h(t)\|^2_{\dbL^2}+ \Vert U_h(t) \Vert^2_{\dbL^2} \, {\rm d}t\Bigr]
+  \frac {\a}{2}\me \bigl[\|X_h(T)\|^2_{\dbL^2}\bigr] \,,
\end{equation}
subject to the equation
\bel{w1013e1a}   
\lt\{
\begin{array}{ll}
{\rm d}X_h(t)=\big[\D_h X_h(t)+ U_h(t)\big]\, {\rm d}t+  \big[X_h(t)+\cR_h\si(t)\big] {\rm d}W(t) \qq &\forall\, t \in [0,T]\,,\\
X_h(0)= \cR_h X_0\,.
\end{array}
\rt.
\ee
Note the use of the Ritz projection, instead of the ${\mathbb L}^2$-projection, which will be useful in Section \ref{stab-1}.
The existence of a unique optimal pair $(X_h^*, U^*_h)$ follows from \cite{Yong-Zhou99}, 
thanks to its characterization via Pontryagin's maximum principle, {\em i.e.},
\begin{equation}\label{pontr1a}  
0 = U^*_h(t) - Y_h(t) \qquad \forall\, t \in (0,T)\, ,
\end{equation}
where the adjoint $(Y_h, Z_h) \in L^2_{{\mathbb F}}\bigl( \Omega; C([0,T]; {\mathbb V}_h)\bigr) \times L^2_{{\mathbb F}}\bigl( 0,T; {\mathbb V}_h\bigr)$ solves the {\bf BSPDE}$_h$
\bel{bshe1a}
\lt\{
\begin{array}{ll}
{\rm d}Y_h(t)= \big[-\D_h Y_h(t)- Z_h(t)+ X^*_h(t)  \big]{\rm d}t
 +Z_h(t) {\rm d}W(t)  \q & \forall\, t \in [0,T]\, ,\\
Y_h(T)= -\a X^*_h(T) \,.
\end{array}
\rt.
\ee
In \cite{Dunst-Prohl16},  error estimates have been obtained for 
$(X_h^*,Y_h, Z_h)$ with the help of a fixed point argument --- which crucially exploits $T>0$ to be {\em sufficiently small}. One main goal in this work is to derive
corresponding estimates for $(X_h^*, Y_h, Z_h,U_h^*)$ for {\em arbitrary} $T>0$ via a variational argument which 
exploits properties of the cost functional $\cJ$ in \rf{w1003e2h}: 
once an estimate for $\me \bigl[\int_0^T\Vert U^*(s) - U^*_h(s)\Vert^2_{{\mathbb L}^2}\, {\rm d}s\bigr]$ stands, 
we use the results from Sections \ref{opt-1} and \ref{opt-2} to derive estimates for
the remaining processes in  $(X_h^*, Y_h, Z_h,U_h^*)$.
%
\bt{rate1}
Under Assumption {\rm(A)}, let $(X^*, U^*)$ be the solution to problem {\bf SLQ}, and $(X^*_h, U^*_h)$ solve 
problem {\bf SLQ$_h$}. 
Then, there exists $C \equiv C(X_0,\si, T)>0 $ independent of $h>0$ such that
\begin{eqnarray*}
{\rm (i)} &&  \me \Bigl[\int_0^T \|U^*(t)-U_h^*(t)\|_{\dbL^2}^2 \rd t\Bigr]  \leq C h^2\, , \\ 
{\rm (ii)} &&  \sup_{0 \leq t \leq T} {\mathbb E} \bigl[\|X^*(t)-X_h^*(t)\|_{\dbL^2}^2\bigr]
       + \me\Bigl[\int_0^T \|X^*(t)-X_h^*(t)\|_{\dbH_0^1}^2 \rd t\Bigr] \leq Ch^2\, ,\\
{\rm (iii)} &&
  \sup_{0 \leq t \leq T} \me \bigl[ \|Y(t)-Y_h(t)\|_{\dbL^2}^2 \bigr]+  \int_0^T {\mathbb E}\Bigl[ \|Y(t)-Y_h(t)\|_{\dbH_0^1}^2+\|Z(t)-Z_h(t)\|_{\dbL^2}^2 \Bigr] \rd t \leq Ch^2\, .
\end{eqnarray*}
\et
We postpone its proof to Section \ref{rate-s}. --- 
In a second step, we propose a temporal discretization of problem {\bf SLQ$_h$} which will be analyzed in 
Section \ref{rate-2}. For this purpose, we use
a mesh $I_\tau$ covering $[0,T]$, and consider step size processes $(X_{h\tau}, U_{h\tau}) \in {\mathbb X}_{h\t} \times {\mathbb U}_{h\t}\subset L^2_\dbF\bigl(0,T;\dbV_h\bigr)\times L^2_\dbF\bigl(0,T;\dbV_h\bigr)$, where
\beq
 {\mathbb X}_{h\t} &:= \lt\{X \in L^2_\dbF (0,T;\dbV_h): \, X(t)=X(t_n), \,\,   \forall t\in [t_n, t_{n+1}),  \,\, n=0,1,\cds, \, N-1\rt\}\, ,\\
{\mathbb U}_{h\t} &:= \lt\{U \in L^2_\dbF (0,T;\dbV_h):\, U(t)=U(t_n),  \,\,    \forall t\in [t_n, t_{n+1}),  \,\, n=0,1,\cds, \, N-1\rt\}\, ,
\eeq
and define for any $X\in \dbX_{h\t}$  and $U\in \dbU_{h\t}$,
\beq
\|X\|_{\dbX_{h\t}}:=\Big(\t \sum_{n=1}^N\me\big[\|X(t_n)\|_{\dbL^2}^2\big]\Big)^{1/2}\,,\quad
\mbox{and} \quad \|U\|_{\dbU_{h\t}}:=\Big(\t \sum_{n=0}^{N-1}\me \big[\|U(t_n)\|_{\dbL^2}^2\big]\Big)^{1/2}.
\eeq
Note that the norms of $\dbX_{h\t}, \,\dbU_{h\t}$ differ: the reason is that for a control system, we should control  `from now on' --- when the current state is known, and the future state is what we care about.

Problem {\bf SLQ$_{h\t}$} then reads as follows: Find an optimal pair 
$(X_{h\t}^*,U_{h\t}^*)\in {\mathbb X}_{h\t} \times {\mathbb U}_{h\t}$ which minimizes the quadratic  cost functional 
\bel{w1003e8}
\bal
\cJ_{\t}(X_{h\t}, U_{h\t})
=\frac {1} 2 \big[\| X_{h\t}\|^2_{\dbX_{h\t}}
+ \Vert U_{h\t} \Vert_{\dbU_{h\t}}^2\big]
+\frac \a 2 \me \bigl[\|X_{h\t}(T)\|_{\dbL^2}^2\bigr]\, ,
\eal
\ee
subject to a forward difference equation
\bel{w1003e7}
\lt\{
\bal
& X_{h\t}(t_{n+1})-X_{h\t}(t_n)= \tau \big[\D_h X_{h\t}(t_{n+1})
 +U_{h\t}(t_n) \big]+ \big[X_{h\t}(t_n)+\cR_h \si(t_n)\big] \D_{n+1}W\\
& \qq\qq\qq\qq\qq\qq\qq\qq\qq\qq \q n=0,1,\cds,N-1\, ,\\
& X_{h\t}(0)= \cR_h X_0\,.
\eal
\rt.
\ee
The following result states a {(discrete)} Pontryagin-type maximum principle for the uniquely solvable problem {\bf SLQ$_{h\tau}$}. We mention the appearing
mapping $K_{h\tau}$ in \rf{w1212e3}$_2$, which is used to give the discrete optimality condition \rf{w1003e12}; see also Remark \ref{w829r1}. 
 In the sequel, this principle is
 used to verify
rates of convergence for the solution to problem {\bf SLQ$_{h\t}$} towards
the solution to  {\bf SLQ}$_h$.

\bt{MP}
Let $\ds A_0 :=\lt(\mathds{1}-\tau \D_h\rt)^{-1}$.
The unique optimal pair $(X^*_{h\t}, U^*_{h\t})\in {\mathbb X}_{h\tau} \times {\mathbb U}_{h\tau}$ of
problem {\bf SLQ$_{h\t}$}  solves the following equalities for $n=0,1,\cds,N-1$:
\bel{w1212e3}
\lt\{
\bal
X_{h\t}^*(t_{n+1})=&A_0^{n+1}\prod_{j=1}^{n+1}\lt(1+ \D_j W\rt)X_{h\t}({0})+\t \sum_{j=0}^{n}A_0^{n+1-j}\prod_{k=j+2}^{n+1}\lt(1+ \D_kW\rt) U_{h\t}^*(t_j)\\
  &+\sum_{j=0}^n A_0^{n+1-j}\prod_{k=j+2}^{n+1}\lt(1+\D_kW \rt)\cR_h\si(t_j)\D_{j+1}W\, , \\ 
\lt( K_{h\t}X^*_{h\t}\rt)(t_n):=&-\t \me\Big[ \sum_{j=n+1}^N A_0^{j-n} \prod_{k={n+2}}^j(1+\D_kW) X^*_{h\t}(t_j)  \Big|\mf_{t_n} \Big]  \\
&-\a \me \Big[ A_0^{N-n} \prod_{k=n+2}^N(1+\D_kW)  X^*_{h\t}(T) \Big|\mf_{t_n}\Big]\,,\\
X^*_{h\t}(0)=&\cR_h X_0\, ,
\eal
\rt.
\ee
 together with the discrete optimality condition
\bel{w1003e12}     
U^*_{h\t}(t_n)-\lt( K_{h\t} X^*_{h\t}\rt)(t_n)=0 \, .
\ee
\et
We verify this characterization of the unique minimizer of problem {\bf SLQ$_{h\t}$} in Section \ref{rate-2}, and then 
use it to prove our second main result.
\bt{rate2}
Under Assumption {\rm(A)}, let $(X_h^*,U^*_h)$ be the solution to problem {\bf SLQ$_h$}, and $(X^*_{h\t}, U^*_{h\t})$ solves problem {\bf SLQ$_{h\t}$}. Then, there exists $C \equiv C(X_0,\si, T)>0$ independent of $h, \tau>0$ such that
\begin{eqnarray*}
{\rm (i)} && \sum_{k=0}^{N-1}\me \Bigl[\int_{t_k}^{t_{k+1}}\lt\|U^*_h(t)-U^*_{h\t}(t_k)\rt\|_{\dbL^2}^2\, {\rm d} t \Bigr]
     \leq C \t \, ; \\ 
{\rm (ii)}&& \max_{0 \leq k \leq N}  {\mathbb E} \bigl[\|X_h^*(t_k)-X_{h\t}^*(t_k)\|_{\dbL^2}^2 \bigr]
+\tau \sum_{k=1}^N {\mathbb E}  \bigl[\|X_h^*(t_k)-X^*_{h\t}(t_k)\|_{\dbH_0^1}^2\bigr]
 \leq C\t  \, .
\end{eqnarray*}

\br{w829r1} 
{\bf 1.}~The strategy of proof in Section \ref{rate-2} for Theorems \ref{MP} and \ref{rate2} can also be used to obtain the results in \cite{Prohl-Wang20}, where the noise in {\bf SPDE}$_h$ \rf{sde} --- as part of the optimality conditions for {\bf SLQ}$_{h}$ ---  is  {\rm additive} instead; but not vice versa. In the setting of \cite{Prohl-Wang20}, the coupled {\bf FBSPDE}$_h$ \rf{w1013e1a}--\rf{bshe1a} is replaced by 
\bel{w902e3}
\lt\{
\bal
&{\rm d}X_h(t)=\bigl[\Delta_h X_h(t)+ Y_h(t)\bigr]\, {\rm d}t+ {\mathcal R}_h \si(t) {\rm d}W(t) \q \forall\, t \in [0,T]\,,\\
&{\rm d}Y_h(t)= \bigl[-\Delta_h Y_h(t)+ X_h(t) \bigr]{\rm d}t+Z_h(t) {\rm d}W(t)  \q  \forall\, t \in [0,T]\, ,\\
&X_h(0)= \Pi_h X_0\, ,\qquad Y_h(T)=-\a X_h(T)\, ,
\eal
\rt.
\ee
and to solve the corresponding problem to {\bf SLQ}$_{h\t}$ is equivalent to  solving the following forward-backward 
stochastic difference equation for
$ n =0,1,\cds,N-1$,
\bel{w1212e30}
\lt\{
\bal
& [\mathds{1} - \tau \Delta_h]X_{h\t}(t_{n+1})= X_{h\t}(t_n)+ \tau Y_{h\t}(t_n)+{\mathcal R}_h \si(t_n) \D_{n+1}W\, ,\\
&[\mathds{1} - \tau \Delta_h]Y_{h\t}(t_n) = {\mathbb E}\big[ Y_{h\t}(t_{n+1})- {\tau} X_{h\t}(t_{n+1})\bigl\vert
{\mathcal F}_{t_n}\big]\, , \\
&X_{h\t}(0)=\Pi_h X_0\, ,\qquad Y_{h\t}(T)=-\a X_{h\t}(T)\, .
\eal
\rt.
\ee
based on which the optimal strong error estimates in \cite[Theorem 4.2]{Prohl-Wang20} are shown.

{
{\bf 2.}~Equation
\eqref{w1212e3}$_1$ is the time-implicit approximation of {\bf SPDE$_h$} \eqref{w212e3}$_1$, while $K_{h\t}X^*_{h\t}$ is 
an approximation of $Y_h$ to {\bf BSPDE}$_h$ \eqref{w212e3}$_2$. In fact, \eqref{w1212e3}$_2$ is different from the 
temporal discretization of \eqref{w212e3}$_2$ via the implicit Euler method; see
also Lemma \ref{w229l1} for $Y_h$'s approximation based on implicit Euler method. 
In Lemma \ref{w301l1}, we {estimate the difference} of these two approximations.
When {\bf SPDE} \rf{w1013e1} is driven by {\em additive} noise ({\em i.e.}, $\si(t)\rd W(t)$), $K_{h\t}X^*_{h\t}$ defined in 
\rf{w1212e3}$_2$ is just $Y_{h\t}$ in \rf{w1212e30} by changing $\Pi_hX_0$ to $\cR_hX_0$,  and Theorem \ref{MP} turns
to \cite[Theorem 4.2]{Prohl-Wang20}; see also item {\bf 1}.

{\bf 3.}~Compared to {\bf BSPDE}$_h$ \rf{w1212e30}$_2$, the adjoint equation \rf{bshe1a}$_2$ for problem {\bf SLQ}$_h$ 
contains $Z_h$ in the drift term. This induces extra difficulties (such as to estimate $\me[\|Z_h(t)-Z_h(s)\|_{\dbL^2}^2]$, $s,t\in[0,T]$) when deducing a 
convergence rate for the temporal discretization \rf{bshe1a}$_2$. In this work,  without any extra assumptions on data, 
we adopt the stochastic Riccati equation to overcome these difficulties; see Section \ref{stab-1} for further details.

}
\er

\et
The optimality system \eqref{w1212e3}--\eqref{w1003e12} is still not amenable to an actual implementation, but serves as a key step towards the practical Algorithm \ref{alg1}, which approximately solves {\bf SLQ}$_{h\t}$; its convergence will be shown in Section \ref{numopt}.

\subsection{{Spatial} semi-discretization {\bf SLQ}$_{h}$: Proof of Theorem \ref{rate1}}\label{rate-s} 

We remark that by \eqref{w212e3}$_1$, $X^*$ may be written as $X^* = {\mathcal S}(U^*)$, 
where 
\beq
{\mathcal S}: L^2_\dbF (0,T; \dbL^2) \rightarrow L^2_{{\mathbb F}}\bigl(\Omega; C([0,T]; {\mathbb H}^1_0)\bigr) \cap
L^2_\dbF(0,T; \dbH_0^1\cap\dbH^2)
\eeq
is the bounded  `control-to-state' map. If $X_0 \equiv 0$ and $\si \equiv 0$, we denote this solution map by $\cS^0$. Moreover, we introduce the  reduced functional
\beq
\widehat{\mathcal J}: L^2_{{\mathbb F}}(0,T; {{\mathbb L}^2}) \rightarrow {\mathbb R} \qquad \mbox{via} \qquad \widehat{\mathcal J}(U ) = {\mathcal J}\bigl(\cS(U), U\bigr)\,,
\eeq
where $\cJ$ is defined in \eqref{w1003e2}.
The solution to equation \eqref{w212e3}$_2$ may be written in the form $(Y,\,Z) = (\cT^1(X^*),\,\cT^2(X^*))$, where 
\beq
\cT^1: L^2_{{\mathbb F}}\big(\Omega; C([0,T]; \dbL^2)\big) \rightarrow L^2_\dbF\big(\Omega; C([0,T]; {\mathbb H}^1_0)\big) \cap L^2_\dbF(0,T; \dbH^1_0 \cap \dbH^2)\, ,
\eeq
\beq
\cT^2: L^2_{{\mathbb F}}\big(\Omega; C([0,T]; \dbL^2)\big) \rightarrow  L^2_\dbF(0,T; {\mathbb H}^1_0)\, ,
\eeq
which are both bounded.

%
\bl{w212l1}
For every $U \in L^2_\dbF (0,T; {\dbL^2}) $, the Fr\'echet derivative
$D\h \cJ(U)$  is a bounded operator on $L^2_{{\mathbb F}} (0,T; {\dbL^2})$
which takes the form
\begin{equation}\label{derivative-cont}
D \widehat{\mathcal J}(U) = U -\cT^1\bigl(\cS(U)\bigr) \, .
\end{equation}
\el

\begin{proof}
By \eqref{stochheat1}, the stability of {\bf SPDE} \eqref{w212e3}$_1$, we have
\beq
\me \bigl[\| \cS^0(V)(t)\|^2_{\dbL^2}\bigr]\leq C \|V\|^2_{L^2_\dbF (0,T;\dbL^2)} \qq \forall~t\in [0,T],\, V\in L^2_\dbF(0,T;\dbL^2)\,.
\eeq
Hence, for any $U,\,V\in L^2_\dbF(0,T;\dbL^2)$, applying the fact $\cS(U+V)=\cS(U)+\cS^0(V)$, we can get
\beq
\bal
&\h\cJ(U+V)-\h\cJ(U)\\
 &\q-\Big[(\cS(U),\cS^0(V))_{L^2_\dbF(0,T;\dbL^2)}+\bigl((\a \cS(U)(T)),\cS^0(V)(T)\bigr)_{L^2_{\mf_T}(\O;\dbL^2)}+(U,V)_{L^2_\dbF(0,T;\dbL^2)}\Big]\\
&=\frac 1 2 \Big[ \big\| \cS^0(V)\big\|^2_{L^2_\dbF(0,T;\dbL^2)}+\|V\|^2_{L^2_\dbF(0,T;\dbL^2)}
+ \a\big\| \cS^0(V)(T)\big\|^2_{L^2_{\mf_T}(\O;\dbL^2)}\Big].
\eal
\eeq
On the other side, It\^o's formula to $(\cS^0(V),\cT^1(\cS(U)))_{\dbL^2}$ yields to
\beq
\bal
\bigl( \cS(U),\cS^0(V) \bigr)_{L^2_\dbF(0,T;\dbL^2)}+ \bigl( \a\cS(U)(T),\cS^0(V)(T) \bigr)_{L^2_{\mf_T}(\O;\dbL^2)}
=-\bigl(\cT^1\lt(\cS(U)\rt), V\bigr)_{L^2_\dbF(0,T;\dbL^2)}.
\eal
\eeq
Subsequently, by the definition of Fr\'echet derivative and above three inequalities,
we prove the desired
result.
\end{proof}

By the unique solvability property of \eqref{w1013e1a}, we associate to this equation the  bounded solution operator
$${\mathcal S}_h: L^2_{{\mathbb F}}(0,T; \dbV_h) \rightarrow L^2_{{\mathbb F}}\bigl(\Omega; C([0,T]; {\mathbb V}_h)\bigr)\,,$$
which allows to introduce the reduced functional
\begin{equation}\label{red-conv}
\widehat{\cJ}_h: L^2_{{\mathbb F}}(0,T; \dbV_h) \rightarrow {\mathbb R}\,, \qquad \mbox{via} \qquad
\widehat{\cJ}_h (U_h) = \cJ \bigl({\mathcal S}_h(U_h), U_h\bigr)\, ,
\end{equation}
where $\cJ$ is defined in \eqref{w1003e2}.
The solution pair to equation \eqref{bshe1a} may be written as $(Y_h,\,Z_h) = \lt(\cT^1_h(X^*_h),\,\cT^2_h(X^*_h)\rt)$, where 
\beq
\cT^1_h: L^2_{{\mathbb F}}\big(\Omega; C([0,T]; {\mathbb V}_h) \big) \rightarrow L^2_{{\mathbb F}}\big(\Omega; C([0,T]; {\mathbb V}_h) \big)\, ,\q
\cT^2_h: L^2_{{\mathbb F}}\big(\Omega; C([0,T]; {\mathbb V}_h) \big) \rightarrow L^2_{{\mathbb F}}(0,T; {\mathbb V}_h )\, .
\eeq

\ms
We are now in a position to prove Theorem \ref{rate1}.

\begin{proof}[\bf {Proof of Theorem \ref{rate1}}] 

For every $U_h \in L^2_{{\mathbb F}}(0,T; {{\mathbb V}_h} ) $, the Fr\'echet derivative 
${D\h\cJ_h(U_h)}$ is a bounded operator (uniformly in $h$) on $L^2_{{\mathbb F}}(0,T; \dbV_h) $, and has the form
\begin{equation}\label{derivative-semidisc}
D {\widehat\cJ}_h (U_h)=U_h- {\mathcal T}_h\bigl({\mathcal S}_h (U_h)\bigr)\, ,
\end{equation}
which can be deduced by the similar procedure as that in Lemma \ref{w212l1}.
Let $U_h \in L^2_{{\mathbb F}}(0,T; {{\mathbb V}_h})$ be arbitrary; it is due to the quadratic structure of the reduced functional \eqref{red-conv} that
$$ \bigl( D^2 \widehat{\mathcal J}_h(U_h) R_h,R_h \bigr)_{L^2_\dbF(0,T;\dbL^2)} 
\geq \Vert R_h\Vert^2_{L^2_\dbF(0,T; {\mathbb L}^2)}
\quad \forall\, R_h \in 
L^2_\dbF(0,T; \dbV_h )\, .$$
As a consequence, on putting $R_h = U^*_h - \Pi_h U^*$,
\beq
\bal
& \Vert U^*_h - \Pi_h U^*\Vert_{L^2_\dbF(0,T; {\mathbb L}^2)}^2
\leq  \big( D^2 \widehat{\mathcal J}_h(U_h) (U^*_h - \Pi_h U^*),U^*_h - \Pi_h U^*\big)_{L^2_\dbF(0,T;\dbL^2)}\\ 
&\qquad =  \bigl( D\widehat{\mathcal J}_h(U^*_h) ,U^*_h - \Pi_h U^*\bigr)_{L^2_\dbF(0,T;\dbL^2)} - 
\bigl( D\widehat{\mathcal J}_h(\Pi_h U^*) ,U^*_h - \Pi_h U^*\bigr)_{L^2_\dbF(0,T;\dbL^2)}\, .
\eal
\eeq
Note that $D\widehat{\mathcal J}_h(U^*_h) = 0$ by \eqref{pontr1a}, as well as $D \widehat{\mathcal J}(U^*) = 0$ by \eqref{w212e4}, such that the last line equals
\beq
\bal
=& \Bigl[ \big( D\widehat{\mathcal J}(U^*) ,U^*_h - \Pi_h U^*
\big)_{L^2_\dbF(0,T;\dbL^2)} - \big( D \widehat{\mathcal J}
(\Pi_h U^*),U^*_h - \Pi_h U^* \big)_{L^2_\dbF(0,T;\dbL^2)}\Bigl] \\ 
&+  \Bigl[\big(  D \widehat{\mathcal J} (\Pi_h U^*), U^*_h - \Pi_h U^*\big)_{L^2_\dbF(0,T;\dbL^2)} 
- \big( D\widehat{\mathcal J}_h(\Pi_h U^*) ,U^*_h - \Pi_h U^*\big)_{L^2_\dbF(0,T;\dbL^2)}\Bigr]\,.
\eal
\eeq
Hence, 
\bel{w229e1}
\bal
&  \Vert U^*_h - \Pi_h U^*\Vert_{L^2_\dbF(0,T; \dbL^2)}^2\\
&\quad \leq2\Big(  \| D\h\cJ(U^*)-D\h\cJ(\Pi_hU^*) \|^2_{L^2_\dbF(0,T;\dbL^2)} 
       + \|  D\h\cJ(\Pi_hU^*)-D\h\cJ_h(\Pi_hU^*) \|^2_{L^2_\dbF(0,T;\dbL^2)}  \Big)\\
&\quad =:2(I+II)\,.       
\eal
\ee
We use \eqref{derivative-cont} to bound $I$ as follows,
\beq
I \leq &2 \Bigl(\lt\|U^* - \Pi_h U^*\rt\|^2_{L^2_\dbF(0,T;\dbL^2)}
       + \lt\|\cT^1\big(\cS(\Pi_h U^*) \big) -\cT^1\big(\cS(U^*) \big)\rt\|^2_{L^2_\dbF(0,T;\dbL^2)}\Bigr) \,.
\eeq
By  stability properties (see also \eqref{vari-2}) for {\bf BSPDE} \eqref{w212e3}$_2$, as well as {\bf SPDE} \eqref{w212e3}$_1$ 
(see also \eqref{stochheat1}), the last term in the above inequality reads
\bel{estima}
\bal
 \leq& C \Bigl(\Vert \big(\cS(U^*)-\cS(\Pi_h U^*)\big)(T)\Vert^2_{L^2_{\mf_T}(\O;\dbL^2)} 
            +  \Vert \cS(U^*)-\cS(\Pi_h U^*)\Vert^2_{L^2_\dbF(0,T; \dbL^2)}\Bigr)\\ 
\leq& C \Vert U^* - \Pi_h U^*\Vert^2_{L^2_\dbF(0,T; \dbL^2)}\, .
\eal
\ee
By optimality condition \eqref{w212e4},  and the regularity properties of the solution
of {\bf FBSPDE} \eqref{w212e3}, we know that already
$U^* \in L^2_{\mathbb F}(0,T; {\mathbb H}^1_0)$; as a consequence, the right-hand side of \eqref{estima} is bounded by $Ch^2$.

We use the representation \eqref{derivative-semidisc} to bound $II$ via
\bel{w229e2}
\bal
II \leq& \lt\|\cT^1\big(\cS(\Pi_h U^*) \big) - \cT^1_h\big( \cS_h(\Pi_h U^*)\big)\rt\|^2_{L^2(0,T;\dbL^2)} \\
\leq & 2 \big[\lt\|\cT^1\big(\cS(\Pi_h U^*) \big) - \cT^1 \big( \cS_h(\Pi_h U^*)\big)\rt\|^2_{L^2_\dbF(0,T;\dbL^2)}\\
 &\qq+ \lt\|\cT^1\big(\cS_h(\Pi_h U^*) \big) - \cT^1_h\big( \cS_h(\Pi_h U^*)\big)\rt\|^2_{L^2_\dbF(0,T;\dbL^2)}\big]\\
=:& 2 (II_1+II_2)\,.
\eal
\ee
In order to bound $II_{1}$, we use stability properties for {\bf BSPDE} \eqref{w212e3}$_2$, in combination with the error estimate \eqref{esti-space1} for \eqref{sde} to conclude
$$II_{1} \leq C \Bigl( \big\Vert \cS(\Pi_h U^*) - \cS_h(\Pi_h U^*) \big\Vert^2_{L^2_\dbF(0,T; \dbL^2)}
+\Vert \cS(\Pi_h U^*)(T) - \cS_h(\Pi_h U^*)(T) \big\Vert^2_{L^2_{\mf_T}(\O;\dbL^2)}\Bigr) \leq Ch^2\, . $$
In order to bound $II_{2}$, we use the error estimate for {\bf BSPDE} \eqref{w212e3}$_2$, and estimate for {\bf SPDE}$_h$
\rf{w1013e1} with $U_h=\Pi_h U^*$ to find
\beq
\bal
II_{2} 
&\leq C h^2\Bigl(\|\cT^1(\cS_h(\Pi_hU^*))\|^2_{L^\infty_\dbF(0,T;L^2(\O;\dbH_0^1))\cap L^2_\dbF(0,T;\dbH_0^1\cap\dbH^2)}
                 +\|\cT^2(\cS_h(\Pi_hU^*))\|^2_{L^2_\dbF(0,T;\dbH_0^1)}\Bigr)\\
&\leq C h^2\Bigl(\|\cS_h(\Pi_hU^*)(T)\|^2_{L^2_{\mf_T}(\O;\dbH_0^1)}
                 +\|\cS_h(\Pi_hU^*)\|^2_{L^2_\dbF(0,T;\dbL^2)}\Bigr)\\
&\leq C h^2\Bigl(\|X_0\|_{\dbH_0^1}^2+\|Y\|^2_{L^2_\dbF(0,T;\dbL^2)}
                 +\|\si \|^2_{L^2_\dbF(0,T;\dbH_0^1)}\Bigr)\\ 
&\leq Ch^2\, .
\eal
\eeq
We now insert these estimates into \eqref{w229e1}, and utilize the optimal condition \eqref{w212e4} to obtain the bound
\beq
\|U^*-U^*_h\|^2_{L^2_\dbF(0,T;\dbL^2)} 
\leq 2\Bigl( \|U^*-\Pi_hU^*\|^2_{L^2_\dbF(0,T;\dbL^2)} +\| U^*_h - \Pi_h U^*\|_{L^2_\dbF(0,T; \dbL^2)}^2\Bigr) \leq Ch^2\, . 
\eeq
This just part (i) of the theorem. 

\ms
Since $U^* \in L^2_{\dbF}\big(0,T;\dbH_0^1\big)$, and (i),
the estimates (ii) and (iii) can be deduced as \eqref{esti-space1} and Theorem \ref{w0911t1} respectively.
\end{proof}

\subsection{Some $h$-independent stability bounds {and time regularity} properties for the solutions to {\bf SPDE}$_h$ and {\bf BSPDE}$_h$}\label{stab-1}

The following lemmata validate $h$-independent bounds and time regularity in relevant norms for the solution
$X^*_h$ of {\bf SPDE}$_h$ \eqref{w1013e1a} where $U_h \equiv U^*_h$,  and
of $(Y_h,\,Z_h)$ which solves {\bf BSPDE$_h$} \eqref{bshe1a}, which are both crucial to derive the convergence rate
for proposed discretization.
To obtain time regularity, two types of assumptions are usually made on terminal conditions: (1) Malliavin differentiability 
(see e.g.~\cite{Hu-Nualart-Song11,Wang16}), or (2) Markovianity (see e.g.~\cite{ZhangJF04,Wang20}).
In this work, due to the optimal control framework,
we adopt the state feedback strategy to show $h$-independent stability bounds of the optimal pair $(X^*_h,U^*_h)$ to problem {\bf SLQ}$_h$ and $(Y_h,Z_h)$ as well as time regularity of $Z_h$. 
Specifically, we introduce the stochastic Riccati equation \rf{riccati} in combination with a backward ODE \rf{ode} and, in particular, prove that the (stochastic) Riccati operator $P_h$
that solves \rf{riccati} may be bounded uniformly in $h$; see Lemma \ref{w229l4}.
Then, we apply the state feedback control $U^*_h=-P_hX^*_h-\f_h$ to deduce $h$-independent stability bounds 
for the tuple $X_h^*$, then for $U^*_h$ and $(Y_h,Z_h)$. Also by the aforementioned feed back control, we further conclude 
the Malliavin differentiability of $X^*_h, U^*_h, Y_h$, obtain $h$-independent bounds for it, and thus conclude the needed time
regularity of $Z_h$; see (\ref{z19}).


\ms

To begin with, we introduce a family of {\bf SLQ}$_h$ problems, parametrized by $t \in [0,T]$; for this purpose, we consider the controlled {\bf SPDE}$_h$ 
\bel{see}
\lt\{
\bal
&\rd X_h(s)=\big[\D_h X_h(s)+ U_h(s)\big]\rd s+  \big[X_h(s)+\cR_h\si(s)\big] \rd W(s) 
     \qq \forall\, s \in [t,T]\,,\\
&X_h(t)= \cR_h X_t
\eal
\rt.
\ee
with {$X_t\in \dbH_0^1$}, and the (parametrized) cost functional
\bel{cost}
\cJ_h(t,X_t;U_h):=\frac 1 2 \me \Bigl[\int_t^T \|X_h(s)\|^2_{\dbL^2}+\|U_h(s)\|^2_{\dbL^2} \rd s \Bigr]
    +\frac \a 2 \me \bigl[\|X_h(T)\|^2_{\dbL^2}\bigr]\,.
\ee
We define the value function as follows:
\bel{w229e11}
V_h(t,X_t):=\inf_{U_h\in L^2_{\dbF}(t,T;\dbL^2)}\cJ_h(t,X_t;U_h)\,.
\ee
Obviously   $\cJ(X_h,U_h)=\cJ_h(0,X_0;U_h)$, and $\cJ(X^*_h,U^*_h)=V_h(0,X_0)$.
The stochastic Riccati equation related to {\bf SLQ$_h$} then reads:
\bel{riccati}
\lt\{
\bal
&P_h'(t)+P_h(t) \D_h + \D_h P_h(t)+ P_h(t)+\mathds{1}_h
  - P_h(t)P_h(t)=0\q \forall~ t\in [0,T]\,,\\
&P_h(T)=\a\mathds{1}_h \,,
\eal
\rt.
\ee
and we consider a backward ODE, 
\bel{ode}
\lt\{
\bal
&\f_h'(t)+\big[\D_h-P_h(t)\big]\f_h(t)+  P_h(t)\cR_h\si(t)=0 \q \forall~ t\in [0,T]\,,\\
&\f_h(T)=0 \,.
\eal
\rt.
\ee
Here $\mathds{1}_h$ denotes the identity operator on $\dbV_h$.
By \cite[Chapter 6, Theorems 6.1 \& 7.2]{Yong-Zhou99}, we know that the stochastic Riccati equation \eqref{riccati} 
admits a unique solution
$P_h\in C\bigl([0,T];\cL(\dbV_h;\dbV_h)\bigr)$  which is nonnegative, symmetric, subsequently \rf{ode} has a unique solution $\f\in C([0,T];\dbV_h)$, and that
\bel{w229e12}
\bal
V_h(t,X_t)=&\frac 1 2 \bigl(P_h(t)\cR_hX_t, \cR_hX_t \bigr)_{\dbL^2}+  \bigl(\f_h(t),\cR_hX_t \bigr)_{\dbL^2}\\
&+\frac 1 2 \int_t^T\big[ \bigl(P_h(s)\cR_h\si(s),\cR_h \si(s) \bigr)_{\dbL^2}+\|\f_h(s)\|^2_{\dbL^2}\big] \rd s\, .
\eal
\ee

\bl{w229l4}
Let $P_h$ be the solution of \eqref{riccati}, and $\f_h$ solves \rf{ode}. Then there exists a constant $C>0$ 
independent of $h>0$ such that
\begin{eqnarray*}
{\rm (i)}&&
\ds \sup_{t\in[0,T]}\|P_h(t)\|_{\cL\big(\dbL^2|_{\dbV_h};\dbL^2|_{\dbV_h}\big)}\leq C\,, \\
{\rm (ii)}&& \ds \sup_{t\in [0,T]}\|\f_h(t)\|_{\dbL^2}^2+ \int_0^T  \|\nb \f_h(s)\|^2_{\dbL^2}+\bigl(P_h(s)\f_h(s),\f_h(s)\bigr)_{\dbL^2}  \rd s \leq C\|\si\|_{L^2_\dbF(0,T;\dbH_0^1)}^2 \,.
\end{eqnarray*}

\el

\begin{proof}
{\bf 1)} We consider problem {\bf SLQ}$_h$
for \rf{see} with $\si \equiv 0$.
In this case, we denote the solution to \rf{see} by $X_h^0 \equiv X_h^0(\cd;\cR_h X_t, U_h)$.
For 
$U_h\equiv 0$,  It\^o's formula then yields
\beq
\sup_{s\in[t,T]}\me \bigl[\|X_h^0(s;\cR_hX_t, 0)\|_{\dbL^2}^2\bigr] \leq e^{ T}\me\bigl[\|\cR_hX_t\|_{\dbL^2}^2\bigr].
\eeq
Hence, \eqref{w229e11} and \eqref{w229e12} lead to
\beq
\bigl(P_h(t)\cR_hX_t, \cR_hX_t \bigr)_{\dbL^2}= 2V_h(t,X_t)\leq  C\me\bigl[\|\cR_hX_t\|_{\dbL^2}^2\bigr] \, ,
\eeq
which, together with the facts that $P_h$ is nonnegative and $\cR_h$ is surjective, implies assertion (i).

{\bf 2)} To verify (ii), we infer from \rf{ode} that
\begin{eqnarray*}
&&\|\f_h(t)\|_{\dbL^2}^2+2\int_t^T\lt[ \|\nb \f_h(s)\|^2_{\dbL^2}+\lt(P_h(s)\f_h(s),\f_h(s)\rt)_{\dbL^2} \rt]\rd s\\
&&\qquad \leq \int_t^T \|\f_h(s)\|_{\dbL^2}^2\rd s+\int_t^T \|P_h(s)\cR_h\si(s)\|_{\dbL^2}^2\rd s\,.
\end{eqnarray*}
Then Gronwall's inequality and (i) settle the assertion. 
\end{proof}

The following lemma collects bounds for the solution $X^*_h$ 
of {\bf SPDE}$_h$ \eqref{w1013e1a} with $U_h \equiv U^*_h$, by exploiting the state
feedback representation \rf{w1011e2} of the optimal control $U^*_h$, and the bounds for the stochastic 
Riccati operator in Lemma \ref{w229l4}, in particular.

\bl{w229l2}
Let $X^*_h$ solve {\bf SPDE}$_h$ \eqref{w1013e1a} with $U_h \equiv U^*_h$. Then for any $t \in [0,T]$,
 $X^*_h(t)\in \dbD^{2,2}(\dbL^2)$, and {there exists an $h$-independent constant $C>0$}  such that
{\small \bel{w826e1}
\lt\{
\bal
&\sup_{t\in[0,T]}\me \bigl[\|X^*_h(t)\|_{\dbL^2}^2\bigr]+\sup_{\th\in[0,T]}\sup_{t\in[\th,T]}\me \bigl[\|D_\th X^*_h(t)\|_{\dbL^2}^2\bigr]
+\sup_{\th,\mu \in[0,T]}\sup_{t\in[\mu\vee\th,T]}\me \bigl[\|D_\mu D_\th X^*_h(t)\|_{\dbL^2}^2\bigr]\\
&\qq\qq\qq\qq\qq\qq\qq\qq\qq\qq\leq C \lt(\|X_0\|_{\dbH_0^1}^2+\|\si\|_{C([0,T];\dbH_0^1)}^2\rt)\,,\\
&\sup_{t\in [0,T]}\me \bigl[\|\nb X^*_h(t)\|_{\dbL^2}^2\bigr] +\me\Bigl[\int_0^T\|\D_h X^*_h(t)\|_{\dbL^2}^2\rd t\Bigr]
\leq C \lt(\|X_0\|_{\dbH_0^1}^2+\|\si\|_{C([0,T];\dbH_0^1)}^2\rt)\,,\\
&\sup_{\th\in[0,T]}\sup_{t\in [\th,T]}\me\bigl[\|\nb D_\th X^*_h(t)\|_{\dbL^2}^2 \bigr]+\sup_{\th\in[0,T]}\me\Bigl[\int_\th ^T\|\D_h D_\th X^*_h(t)\|_{\dbL^2}^2\rd t\Bigr]\\
&\qq\qq\qq\qq\qq\qq\qq\qq\qq\qq\leq C \lt(\|X_0\|_{\dbH_0^1}^2+\|\si \|_{C([0,T];\dbH_0^1)}^2\rt)\,,\\
&\me \bigl[\|X^*_h(t)-X^*_h(s)\|_{\dbL^2}^2\bigr] \leq C|t-s| \lt(\|X_0\|_{\dbH_0^1}^2+\|\si \|_{C([0,T];\dbH_0^1)}^2\rt) \qq t,s\in [0,T]\,,\\
 &\sup_{t\in[\th_1,T]}\me\bigl[\|(D_{\th_1}-D_{\th_2}) X^*_h(t)\|_{\dbL^2}^2 \bigr]
\leq C|\th_1-\th_2|  \lt(\|X_0\|_{\dbH_0^1}^2+\|\si \|_{C^1([0,T];\dbH_0^1)}^2\rt)
\qq \th_2\leq \th_1,\, \\
&\sup_{t\in[0,T]}\me \bigl[\|U^*_h(t)\|_{\dbL^2}^2\bigr] +\Big(\sup_{t\in[0,T]}\me\bigl[\|U^*_h(t)\|_{\dbL^2}^4\bigr]\Big)^{1/2}
   \leq C \lt(\|X_0\|_{\dbH_0^1}^2+\|\si \|_{C([0,T];\dbH_0^1)}^2\rt)\, .
\eal
\rt.
\ee}
\el

\begin{proof}
{\bf 1)} By \cite[Chapter 6, Theorem 6.1, 7.2]{Yong-Zhou99}, the optimal control $U^*_h$ has the following state feedback form
\bel{w1011e2}
U^*_h=-P_h X^*_h-\f_h\, ,
\ee
where $P_h,\f_h$ are solutions to \rf{riccati}, \rf{ode}.
 It\^o's formula for $\|X^*_h\|^2_{\dbL^2}$, in combination with Lemma \ref{w229l4}, (i), yields 
 \beq
 &\me \bigl[\|X_h^*(t)\|_{\dbL^2}^2\bigr]+2\me \int_0^t \Bigl[ \|\nb X_h^*(s)\|_{\dbL^2}^2+\bigl(P_h(s)X_h^*(s),X_h^*(s)\bigr)_{\dbL^2}\Bigr] \rd s\\
 &\quad \leq \|\cR_h X_0\|_{\dbL^2}^2+\int_0^T \Bigl[ \|\f_h(t)\|_{\dbL^2}^2+2\|\cR_h\si(t)\|_{\dbL^2}^2\Bigr] \rd t
      +3\me \Bigl[\int_0^t \|X_h^*(t)\|_{\dbL^2}^2\rd s\Bigr]\,.
 \eeq
 By Gronwall's inequality, Lemma \ref{w229l4}, (ii), and stability properties of ${\mathcal R}_h$ we conclude that
\beq
&\sup_{t\in[0,T]}\me \bigl[\|X^*_h(t)\|_{\dbL^2}^2\bigr]+\me\Bigl[\int_0^T \|\nb X_h^*(t)\|_{\dbL^2}^2 +\bigl(P_h(t)X_h^*(t),X_h^*(t)\bigr) \rd t\Bigr]\\
&\quad \leq C \Bigl(\|\cR_hX_0\|_{\dbL^2}^2+\int_0^T\Bigl[ \|\f_h(t)\|_{\dbL^2}^2+\|\cR_h\si(t)\|_{\dbL^2}^2\Bigr] \rd t \Bigr) \\
&\q \leq C \Bigl(\|X_0\|_{\dbH_0^1}^2+\int_0^T\|\si(t)\|_{\dbH_0^1}^2 \rd t\Bigr)\,.
\eeq 

{\bf 2)} To estimate the second term on the left-hand side of \eqref{w826e1}$_1$, by noting that $ P_h$ and $\f_h$ are deterministic and
\cite[Theorem 2.2.1]{Nualart06}, the first Malliavin derivative of $X_h^*$ exists.  We take the Malliavin
derivative on both sides of \eqref{w1013e1a} with the sate feedback control \rf{w1011e2}, and then
\bel{w301e8}
\lt\{
\bal
&\rd D_\th X^*_h(t)=\big[\D_h- P_h(t) \big] D_\th X^*_h(t)\rd t
   + D_\th X^*_h(t)\rd W(t) \qq \forall\, t\in [\th,T]\,,\\
&D_\th X^*_h(\th)= X^*_h(\th)+\cR_h\si(\th)\, ,\\
&D_\th X^*_h(t)=0 \qq \forall\, t\in [0,\th)\, .
\eal
\rt.
\ee
 It\^o's formula leads to 
\beq
\sup_{\th\in[0,T]}\sup_{t\in[\th,T]}\me \bigl[\|D_\th X^*_h(t)\|_{\dbL^2}^2\bigr]
\leq C\sup_{\th\in[0,T]} \me \bigl[\|X_h^*(\th)+\cR_h\si(\th)\|_{\dbL^2}^2 \bigr]
\leq C\lt[\|X_0\|_{\dbH_0^1}^2+\|\si\|_{C([0,T];\dbH_0^1)}^2\rt]\, .
\eeq 
Hence, $X^*_h\in \dbD^{1,2}(\dbL^2)$. In a similar vein, we can deduce $X^*_h\in \dbD^{2,2}(\dbL^2)$ and
the remaining part of \eqref{w826e1}$_1$.

{\bf 3)} For \rf{w826e1}$_2$, similarly to {\bf 1)}, we estimate
\beq
&\sup_{t\in [0,T]}\me \bigl[\|\nb X^*_h(t)\|_{\dbL^2}^2\bigr]+\me\Bigl[\int_0^T\|\D_h X^*_h(t)\|_{\dbL^2}^2\rd t\Bigr]\\
&\quad \leq C\Big(\|\nb\cR_h X_0\|_{\dbL^2}^2+\me\Bigl[\int_0^T \|P_h(t)X_h^*(t)\|_{\dbL^2}^2+\|\f_h(t)\|_{\dbL^2}^2\
       +\|\nb\cR_h\si(t)\|_{\dbL^2}^2 \rd t\Bigr] \Big)\\
&\quad\leq C \lt( \|X_0\|_{\dbH_0^1}^2+\|\si\|^2_{L^2(0,T;\dbH_0^1)} \rt) \,,      
\eeq
where Lemma \ref{w229l4}, assertion \rf{w826e1}$_1$, and stability properties of ${\mathcal R}_h$  are 
applied.

{\bf 4)} We use It\^o's formula for the solution of 
\rf{w301e8} to find via Lemma \ref{w229l4}, (i) that
\beq
&\me \bigl[\|\nb D_\th X^*_h(t)\|_{\dbL^2}^2\bigr] +\me\Bigl[\int_\th^t\|\D_h D_\th X^*_h(s)\|_{\dbL^2}^2\rd s\Bigr]\\
&\quad \leq \lt( \me\bigl[\|\nb D_\th X^*(\th)\|_{\dbL^2}^2\bigr]+\me\Bigl[\int_\th^T\|P_h(s)D_\th X_h^*(s)\|_{\dbL^2}^2 \rd s\Bigr]\rt)
       +\me\Bigl[\int_\th^t \|\nb D_\th X_h^*(s)\|_{\dbL^2}^2\rd s\Bigr]\\
&\quad \leq C  \lt( \me \bigl[\|\nb X^*(\th)\|_{\dbL^2}^2\bigr] +\|\nb \cR_h\si(\th)\|_{\dbL^2}^2+\me \Bigl[\int_\th^T\|D_\th X_h^*(s)\|_{\dbL^2}^2 \rd s\Bigr]\rt)
       +\me \Bigl[\int_\th^t \|\nb D_\th X_h^*(s)\|_{\dbL^2}^2\rd s\Bigr] \, .   
\eeq
Then, by Gronwall's inequality as well as \rf{w826e1}$_1$, \eqref{w826e1}$_2$, 
this shows \rf{w826e1}$_3$.

{\bf 5)} The verification of estimate \eqref{w826e1}$_4$ can now be deduced via {\bf SPDE}$_h$ \eqref{w1013e1a} and \rf{w826e1}$_1$, \eqref{w826e1}$_2$.

{\bf 6)} We estimate \eqref{w826e1}$_5$. Using equation \eqref{w301e8}, It\^o's formula, and
\eqref{w826e1}$_1$--\eqref{w826e1}$_3$,  we find that
\beq
&\sup_{t\in[\th_1,T]}\me\bigl[\|(D_{\th_1}-D_{\th_2}) X^*_h(t)\|_{\dbL^2}^2\bigr]\\
&\quad \leq C\Bigl(\me\bigl[\|D_{\th_1}X^*_h(\th_1)-D_{\th_2}X^*_h(\th_2)\|_{\dbL^2}^2\bigr]
+ \me\bigl[\|D_{\th_2}X^*_h(\th_1)-D_{\th_2}X^*_h(\th_2)\|_{\dbL^2}^2\bigr] \Bigr)\\
&\quad \leq C\Bigl(  \me \bigl[\|\big(X^*_h(\th_1)+\cR_h\si(\th_1)\big)-\big(X^*_h(\th_2)\|-\cR_h\si(\th_2)\big)\|_{\dbL^2}^2\bigr]\\
&\qquad+(\th_1-\th_2)\me\Bigl[\int_0^T\|\D_h D_{\th_2}X^*_h(s)\|_{\dbL^2}^2\rd s\Bigr]
  +\me\Bigl[\int_{\th_2}^{\th_1}\|D_{\th_2}X^*_h(s)\|_{\dbL^2}^2\rd s\Bigr]\Bigr)\\
&\quad \leq C|\th_1-\th_2|\Big( \|X_0\|_{\dbH_0^1}^2+(\th_1-\th_2)\int_{\th_2}^{\th_1}\|\si'(t)\|_{\dbH_0^1}^2\rd t+\|\si\|_{C([0,T];\dbH_0^1)}^2\Big) \,,
\eeq
which settles the assertion.

{\bf 7)} By \rf{w1011e2}, Lemma \ref{w229l4}, and \eqref{w826e1}$_1$, we 
easily find the first estimate in \eqref{w826e1}$_6$.
%
Then, the application of It\^o's formula to $\|X^*_h\|_{\dbL^2}^4$  leads to the remaining estimate in \eqref{w826e1}$_6$.
\end{proof}

 We now may use these estimates for  $X^*_h$ 
 to bound the solution
$(Y_h,\,Z_h)$ of {\bf BSPDE$_h$} \eqref{bshe1a} --- in which $X^*_h$ appears as well.

\bl{w229l3}
Suppose that $(Y_h,\,Z_h)$ solves {\bf BSPDE$_h$} \eqref{bshe1a} and $I_\t$ is a uniform time mesh of $[0,T]$. Then
there exists a constant $C>0$ independent of $h>0$ such that 
{\footnotesize
\bel{w301e1}
\lt\{
\bal
&\sup_{t\in [0,T]}\me \bigl[\|Y_h(t)\|_{\dbL^2}^2\bigr]
+\me \Bigl[\int_0^T \|\nb Y_h(t)\|_{\dbL^2}^2+\| Z_h(t)\|_{\dbL^2}^2 \rd t \Bigr]
  \leq C \lt(\|X_0\|_{\dbH_0^1}^2+\|\si\|_{L^2(0,T;\dbH_0^1)}^2\rt),\\
&\sup_{t\in [0,T]}\me \bigl[\|\nb Y_h(t)\|_{\dbL^2}^2\bigr]
+\me\Bigl[\int_0^T \|\D_hY_h(t)\|_{\dbL^2}^2+\|\nb Z_h(t)\|_{\dbL^2}^2  \rd t\Bigr]
\leq  C \lt(\|X_0\|_{\dbH_0^1}^2+\|\si \|_{L^2(0,T;\dbH_0^1)}^2\rt),\\
&\|Y_h-\Pi_\t Y_h\|_{L^2_\dbF(0,T;\dbL^2)}^2\leq C \t \lt(\|X_0\|_{\dbH_0^1}^2+\|\si\|_{L^2(0,T;\dbH_0^1)}^2\rt)\,,\\
&\sup_{\th\in[0,T]} \sup_{t\in [\th, T]}\me\bigl[\|D_\th Y_h(t)\|_{\dbL^2}^2\bigr]
      +\sup_{\th,\mu\in[0,T]}\sup_{t\in [\mu\vee\th, T]}\me\bigl[\|D_\mu D_\th Y_h(t)\|_{\dbL^2}^2\bigr]
      +\sup_{\th\in[0,T]}\me\Bigl[\int_\th^T \|D_\th Z_h(t)\|^2_{\dbL^2}\rd t\Bigr]\\
&\qq\qq\qq\qq \qq\qq\qq\qq      
      \leq C \lt(\|X_0\|_{\dbH_0^1}^2+\|\si \|_{C([0,T];\dbH_0^1)}^2\rt)\,,\\
&\sup_{\th\in[0,T]}\sup_{t\in [\th,T]}\me\bigl[\|\nb D_\th Y_h(t)\|_{\dbL^2}^2\bigr]
+\sup_{\th\in[0,T]}\me\Bigl[\int_\th^T  \|\D_h D_\th Y_h(t)\|_{\dbL^2}^2+\|\nb D_\th Z_h(t)\|_{\dbL^2}^2  \rd t\Bigr]\\
&\qq\qq\qq\qq \qq\qq\qq\qq  \leq  C \lt(\|X_0\|_{\dbH_0^1}^2+\|\si\|_{C([0,T];\dbH_0^1)}^2\rt)\, ,\\      
&\me\bigl[\|Z_h(t)-Z_h(s)\|_{\dbL^2}^2\bigr]\leq C |t-s|  \lt(\|X_0\|_{\dbH_0^1}^2+\|\si\|_{C^1([0,T];\dbH_0^1)}^2\rt)  \qq  s,\,t\in[0,T]\,,      
\eal
\rt.
\ee}
where  the (piecewise constant) operator ${\Pi_\t}: L^2_\dbF\big(\O;C([0,T]; \dbV_h)\big)\rightarrow {\mathbb U}_{h\t}$
is defined by
\bel{w828e1}
\Pi_\t U_h(t):=U_h(t_n)  \qquad \forall\, t \in [t_n,t_{n+1}) \qquad n=0,1,\cds, N-1\, .
\ee
\el

\begin{proof}
{\bf 1)} Applying It\^o's formula for $\|\nb Y_h\|_{\dbL^2}^2$ in {\bf BSPDE$_h$} \eqref{bshe1a}, we see that
\beq
\bal
&\sup_{t\in [0,T]}\me\bigl[\|\nb Y_h(t)\|_{\dbL^2}^2\bigr]
+\me\Bigl[\int_0^T  \|\D_hY_h(t)\|_{\dbL^2}^2+\|\nb Z_h(t)\|_{\dbL^2}^2  \rd t\Bigr]\\
&\quad \leq  C\Bigl( \me\bigl[\lt\|\nb  X^*_h(T)\rt\|_{\dbL^2}^2\bigr]
  +\me\Bigl[\int_0^T \|X^*_h(t)\|_{\dbL^2}^2\rd t\Bigr]\Bigr)\,.
\eal
\eeq
By \eqref{w826e1}$_1$, \eqref{w826e1}$_2$ and above inequality, we can get  \eqref{w301e1}$_2$.  Estimate \eqref{w301e1}$_1$ can be obtained in the same vein.

{\bf 2)} {\bf BSPDE$_h$} \eqref{bshe1a} and  \eqref{w826e1}$_1$, \eqref{w301e1}$_1$, \eqref{w301e1}$_2$ lead to
\bel{w301e2}
\bal
&\|Y_h-\Pi_\t Y_h\|_{L^2_\dbF(0,T;\dbL^2)}^2\\
&= \sum_{k=0}^{N-1}\me\Bigl[\int_{t_k}^{t_{k+1}}\bigl\Vert \int_{t_k}^t -\D_hY_h(s)-Z_h(s)+X^*_h(s)  \rd s+\int_{t_k}^t Z_h(s)\rd W(s) \bigr\Vert^2_{\dbL^2}\rd t\Bigr]\\
&\leq C\t \me \Bigl[ \int_0^T \|\D_h Y_h(t)\|_{\dbL^2}^2+\|Z_h(t)\|_{\dbL^2}^2
   +\|X^*_h(t)\|_{\dbL^2}^2  \rd t\Bigr]\\
&\leq C \lt(\|X_0\|_{\dbH_0^1}^2+\|\si\|_{L^2(0,T;\dbH_0^1)}^2\rt)\,,  
\eal
\ee
which is \eqref{w301e1}$_3$.

\ms
{\bf 3)} Noting that $X^*_h(t)\in \dbD^{2,2}(\dbL^2)$,  by \cite[Proposition 5.3]{ElKaroui-Peng-Quenez97} we can get
\bel{w301e5}
\lt\{
\bal
&\rd D_\th Y_h(t)=\big[-\D_h D_\th Y_h(t)- D_\th Z_h(t)+ D_\th X^*_h(t)\big]\rd t   
   +D_\th Z_h(t)\rd W(t)\q t\in [\th, T]\, ,\\
&D_\th Y_h(T)=-\a D_\th X^*_h(T)\,,\\
&D_\th Y_h(t) =0,\,\,D_\th Z_h(t)=0 \q t\in [0,\th)\,,       
\eal
\rt.
\ee
and
\bel{w301e6}
\lt\{
\bal
&\rd D_\mu D_\th Y_h(t)=\big[-\D_h D_\mu D_\th Y_h(t)- D_\mu D_\th Z_h(t)
 + D_\mu D_\th X^*_h(t)\big]\rd t\\
 &\qq\qq\qq\qq\qq        +D_\mu D_\th Z_h(t)\rd W(t)\q t\in [\mu\vee\th, T],\\
&D_\mu D_\th Y_h(T)=-\a  D_\mu D_\th X^*_h(T)\,.      
\eal
\rt.
\ee
Then, It\^o's formula and \eqref{w826e1}$_1$, \eqref{w826e1}$_3$ lead to \eqref{w301e1}$_4$ and \eqref{w301e1}$_5$.

{\bf 4)} To estimate \eqref{w301e1}$_6$, by applying the fact that $Z_h(\cd)=D_\cd Y_h(\cd), \ae$, we arrive at
\bel{w301e10}
\bal
&\me \bigl[\|Z_h(t)-Z_h(s)\|_{\dbL^2}^2\bigr] =\me\bigl[\|D_tY_h(t)-D_sY_h(s)\|_{\dbL^2}^2\bigr]\\
&\quad \leq 2\me\bigl[\|(D_t-D_s)Y_h(t)\|_{\dbL^2}^2\bigr]+2\me\bigl[\|D_s(Y_h(t)-Y_h(s))\|_{\dbL^2}^2\bigr]\,.
\eal
\ee
For the first term on the right side of \rf{w301e10}, It\^o's formula and \eqref{w826e1}$_5$ lead to
\bel{w1023e1}
\bal
&\me\bigl[\|(D_t-D_s)Y_h(t)\|_{\dbL^2}^2\bigr]\\
&\leq C \Big(\me\bigl[\|(D_t-D_s)X^*_h(T)\|_{\dbL^2}^2\bigr]+ \me\Bigl[\int_t^T\|(D_t-D_s)X^*_h(\th)\|_{\dbL^2}^2\rd \th\Bigr]\Big)\\
&\leq C|t-s| \lt(\|X_0\|_{\dbH_0^1}^2+\|\si\|_{C^1([0,T];\dbH_0^1)}^2\rt)\,.
\eal
\ee
Similar to \rf{w301e2},  by virtue of \eqref{w826e1}$_1$, \eqref{w301e1}$_4$, \eqref{w301e1}$_5$
\bel{w1023e2}
\bal
&\me\bigl[\|D_s(Y_h(t)-Y_h(s))\|_{\dbL^2}^2\bigr] \\
&\qq \leq 
C\Bigl( \me\Bigl[\int_s^t\|D_sD_\th Y_h(\th)\|_{\dbL^2}^2\rd \th\Bigr]\\
&\qq\qq+(t-s)\me\Bigl[ \int_s^T \|\D_h D_sY_h(\th)\|^2_{\dbL^2}+\| D_sZ_h(\th)\|^2_{\dbL^2}
+\|D_sX^*_h(\th)\|^2_{\dbL^2} \rd \th \Bigr]\Bigr)\\
& \qq \leq  C|t-s|  \lt(\|X_0\|_{\dbH_0^1}^2+\|\si\|_{C^1([0,T];\dbH_0^1)}^2\rt)\,.
\eal
\ee
Now \rf{w301e1}$_6$ can be derived by \rf{w301e10}--\rf{w1023e2}.
%
\end{proof}

Based on Lemma \ref{w229l2} and Lemma \ref{w229l3}, we can now sharpen our results to the following one.

\bl{w1011l1}
Let $X^*_h$ solve {\bf SPDE}$_h$ \eqref{w1013e1a} with $U_h=U^*_h$, and 
$(Y_h,\,Z_h)$ be the solution to {\bf BSPDE$_h$} \eqref{bshe1a}. Then 
{\small \bel{w1011e1}
\lt\{
\bal
&\sup_{t\in [0,T]}\me \bigl[\|\D_h X^*_h(t)\|_{\dbL^2}^2\bigr]
+\me\Bigl[\int_0^T \|\nb\D_hX^*_h(t)\|_{\dbL^2}^2\rd t\Bigr]
\leq C  \lt(\|X_0\|_{\dbH_2}^2+\|\si\|^2_{L^2(0,T;\dbH^2)}\rt)\, ,\\
&\sup_{t\in [0,T]}\me\bigl[\|\nb\D_h X^*_h(t)\|_{\dbL^2}^2\bigr]
+\me\Bigl[\int_0^T \|\D_h^2X^*_h(t)\|_{\dbL^2}^2\rd t\Bigr]
\leq C\lt( \|X_0\|_{\dbH^3}^2
+\|\si\|^2_{L^2(0,T;\dbH^3)}
\rt)\, ,\\
&\sup_{t\in [0,T]}\me\bigl[\|\D_h Y_h(t)\|_{\dbL^2}^2\bigr]
+\me \Bigl[\int_0^T  \|\nb\D_hY_h(t)\|_{\dbL^2}^2+\|\D_h Z_h(t)\|_{\dbL^2}^2 \rd t
\Bigr] \\
&\qquad \qquad \qquad \qquad \qquad \qquad \qquad \qquad \qquad
\leq C  \lt(\|X_0\|_{\dbH_2}^2+\|\si\|^2_{L^2(0,T;\dbH^2)}\rt)\,,\\
&\me \bigl[\|\nb (X^*_h(t)-X^*_h(s))\|_{\dbL^2}^2\bigr] 
   \leq C |t-s|\lt( \|X_0\|_{\dbH^3}^2+\|\si\|^2_{C([0,T];\dbH_0^1)}+\|\si\|^2_{L^2(0,T;\dbH^3)} \rt) 
   \, ,
\eal \rt.\ee}
for all $t,s \in [0,T]$.

\el

\begin{proof}
Firstly, by definitions of $\D_h$, $\cR_h$, and the fact that $X_0\in \dbH_0^1$, we have
\beq
(\D_h \cR_h X_0,\phi_h)_{\dbL^2}
=-(\nb X_0,\nb\phi_h)_{\dbL^2}
=(\D X_0,\phi_h)_{\dbL^2}=(\Pi_h\D X_0,\phi_h)_{\dbL^2} \q \forall \phi_h\in \dbV_h\,.
\eeq
Hence, we deduce that
\bel{w1023e3}
\D_h \cR_h X_0=\Pi_h\D X_0\,.
\ee
Based on the state feedback control \rf{w1011e2},  It\^o's formula to $\|\D_h X^*_h\|_{\dbL^2}^2$ and \rf{pontr1a}, lead to
\beq
\bal
&\me \bigl[\|\D_h X^*_h(t)\|_{\dbL^2}^2\bigr]
+\me \Bigl[\int_0^t \|\nb \D_hX^*_h(s)\|_{\dbL^2}^2\rd s\Bigr]\\
\leq &  \lt\|\D_h  \cR_hX_0\rt\|_{\dbL^2}^2
 +C\me\Bigl[ \int_0^T  \|\nb U_h(t)\|_{\dbL^2}^2+\|\D_h \cR_h\si(t)\|_{\dbL^2}^2 \rd t
 \Bigr]+\me\Bigl[\int_0^t\|\D_h X^*_h(s)\|_{\dbL^2}^2\rd s\Bigr]\\
= &  \lt\|\D_h  \cR_hX_0\rt\|_{\dbL^2}^2
  +C\me\Bigl[\int_0^T  \|\nb Y_h(t)\|_{\dbL^2}^2+\|\D_h \cR_h\si(t)\|_{\dbL^2}^2 \rd t\Bigr] +\me\Bigl[\int_0^t\|\D_h X^*_h(s)\|_{\dbL^2}^2\rd s\Bigr]\,,\\
\eal
\eeq
which, together with \rf{w1023e3}, \rf{w301e1}$_1$ and Gronwall's inequality, leads to \rf{w1011e1}$_1$. Applying It\^o's formula to 
$\|\D_h Y_h\|_{\dbL^2}^2$ and \rf{w1011e1}$_1$, we can derive \rf{w1011e1}$_3$. In the same vein, on using \rf{w1011e1}$_3$, \rf{w1023e3}, and the $\dbH^1$-stability of the ${\mathbb L}^2$-projection $\Pi_h$ (see e.g.~\cite{Crouzeix-Thomee87}), as well as the fact
\beq
&\|\nb\D_h X_{h\t}(0)\|_{\dbL^2}^2=  \|\nb\D_h \cR_h X_0\|_{\dbL^2}^2
=\|\nb\Pi_h\D X_0\|_{\dbL^2}^2
\leq C\|\D X_0\|_{\dbH_1}^2\leq C \|X_0\|^2_{\dbH^3}\,,
\eeq
thanks to Assumption {\bf  (A)}, and
\beq
&\me \bigl[\int_0^T  \|\nb \D_h \cR_h\si(t)\|_{\dbL^2}^2  \rd t\bigr]
\leq C \me \bigl[\int_0^T  \|\si(t)\|_{\dbH^3}^2  \rd t \bigr]\,,
\eeq
we can derive that
\beq
\bal
&\sup_{t\in[0,T]}\me\bigl[\|\nb\D_h X^*_h(t)\|_{\dbL^2}^2\bigr]
+\me\Bigl[\int_0^T \| \D_h^2X^*_h(t)\|_{\dbL^2}^2\rd t\Bigr]\\
&\quad \leq  C\Bigl( \lt\|\nb\D_h  \cR_hX_0\rt\|_{\dbL^2}^2
  +\me \Bigl[\int_0^T  \|\D_h Y_h(t)\|_{\dbL^2}^2+\|\nb \D_h \cR_h\si(t)\|_{\dbL^2}^2  \rd t\Bigr]\Bigr)\\
&\quad \leq  C\lt( \|X_0\|_{\dbH^3}^2+\|\si\|^2_{L^2(0,T;\dbH^3)}\rt)\,.  
\eal
\eeq
That is \rf{w1011e1}$_2$.

By {\bf SPDE}$_h$ \eqref{w1013e1a}, and \rf{w1011e1}$_1$, \rf{w826e1}$_2$ in Lemma \ref{w229l2}, \rf{w301e1}$_2$
in Lemma \ref{w229l3} as well as \ref{w1024e1},  for $s,t\in [0,T],\, s\leq t$, we can derive that
\beq
\bal
&\me \bigl[\|\nb (X^*_h(t)-X^*_h(s))\|_{\dbL^2}^2\bigr]\\
&\leq C \int_s^t \me\big[ \|\nb\D_hX_h^*(\t)\|_{\dbL^2}^2+\|\nb X_h^*(\t)\|_{\dbL^2}^2
         +\|\nb Y_h(\t)\|_{\dbL^2}^2+\|\nb\cR_h\si(\t)\|_{\dbL^2}^2 \big]\rd \t\\ 
&\leq C |t-s|\lt( \|X_0\|_{\dbH^3}^2+\|\si\|^2_{C([0,T];\dbH_0^1)}+\|\si\|^2_{L^2(0,T;\dbH^3)} \rt)\,. 
\eal
\eeq
That is  \rf{w1011e1}$_4$.
\end{proof}

%

%

\subsection{{Spatio-temporal} discretization {\bf SLQ}$_{h\tau}$: Proof of Theorems \ref{MP} and \ref{rate2}}\label{rate-2}

\begin{proof}[\bf {Proof of Theorem \ref{MP}}]
The proof is similar to that of \cite[Theorem 4.2]{Prohl-Wang20}, and for the completeness, we provide it here;
it consists of two steps.

\ms

\no{\bf 1)}  Recall the definition of $\ds A_0$ in Theorem \ref{MP}.
We define the bounded operators $\G : \dbV_h\rightarrow \dbX_{h\t}$ and $L: \dbU_{h\t}\rightarrow \dbX_{h\t}$
as follows, 
\bel{rep1}
\bal
& \bigl(\G X_{h\t}({0})\bigr)(t_{n})=A_0^n\prod_{j=1}^n\lt(1+\D_j W\rt)X_{h\t}({0})\,,\\
& (LU_{h\t} )(t_{n})=\t \sum_{j=0}^{n-1}A_0^{n-j}\prod_{k=j+2}^n\lt(1+\D_kW\rt) U_{h\t}(t_j)\qquad \forall\, n=1,2,\cds,N\,,
\eal
\ee
respectively, where $(X_{h\tau}, U_{h\tau})$ is an admissible pair in problem {\bf SLQ}$_{h\tau}$; see \rf{w1003e8}--\rf{w1003e7}. We also need $f(\cd)$, which we define as
\beq
f(t_n)=\sum_{j=0}^{n-1} A_0^{n-j}\prod_{k=j+2}^{n}\lt(1+\D_kW \rt)\cR_h\si(t_j)\D_{j+1}W \qquad \forall\, n=1,2,\cds,N\,,
\eeq
%
and use below the abbreviations
\begin{equation}\label{w1003e14}
\h \G {\cR_h X_0} := \lt(\G\cR_h X_0\rt)(T)\, , 
\qquad \h L{U_{h\tau}} := (LU_{h\tau})(T)\,,
\qq \h f:=f(T) \,.
\end{equation}

By \rf{w1003e7}, we can find that
\bel{w1024e2}
X_{h\t}(t_n)=\bigl(\G X_{h\t}(0)\bigr)(t_n)+(LU_{h\t})(t_n)+f(t_n) \qq n=1,2,\cds,N.
\ee

{\bf Claim}: For any $\xi\in {\mathbb X}_{h\t}$, and any $\eta\in L^2_{\mf_T}(\O;\dbV_h)$,
\bel{w1003e01}
\bal
&(L^*\xi)(t_j)=\t \me\Big[ \sum_{n=j+1}^NA_0^{n-j} \prod_{k={j+2}}^n(1+\D_kW) \xi(t_n) \Big|\mf_{t_j} \Big]\,,\\
& (\h L^*\eta)(t_j) = \me \Big[ A_0^{N-j} \prod_{k=j+2}^N(1+\D_kW) \eta \Big|\mf_{t_j}\Big]\q \q j=0,1,\cds,N-1\, .
\eal
\ee

{\em Proof of Claim:} 
{Let $U_{h\t} \in {\mathbb U}_{h\t}$ be arbitrary}. By the definition of $L$ and the fact 
$A_0=A_0^*$, we can calculate that
\beq
&\tau \sum_{n=1}^N \me \big[\big( (LU_{h\t})(t_n),\xi(t_n)  \big)_{ \dbL^2}\big] \\
&\quad = \tau \sum_{j=0}^{N-1} \me\Big[\Bigl( { U_{h\t}(t_j)},  
        \t \me\Big[ \sum_{n=j+1}^NA_0^{n-j} \prod_{k={j+2}}^n(1+\D_kW) \xi(t_n) \bigl \vert\mf_{t_j} \Big] \Bigr)_{\dbL^2}\Big] \, ,
\eeq
which is the first part of the claim. The remaining part can be deduced similarly.
%

\no{\bf 2)} By \rf{w1024e2}, and \eqref{rep1} together with \eqref{w1003e14}, we can rewrite $\cJ_{\t}(X_{h\t} ,U_{h\t} )$ in \rf{w1003e8} as follows:
\beq
\cJ_{\t}(X_{h\t} ,U_{h\t} )
& =\frac 1 2\Big[ \lt( \G {\cR_h X_0} +LU_{h\t}+f ,   \G {\cR_h X_0}+LU_{h\t}+f  \rt)_{L^2_{{\mathbb F}}(0,T; \dbL^2)} 
 + ( U_{h\t} ,U_{h\t})_{L^2_{{\mathbb F}}(0,T; \dbL^2)}\\
 &\qq
   +\a \big(\h\G {\cR_h X_0} +\h LU_{h\t}+\h f, \h\G { \cR_h X_0}+\h LU_{h\t} +\h f \big)_{L^2_{\mf_T}(\Omega; \dbL^2)}\Big]\, ,
\eeq   
Rearranging terms then further leads to
   \beq
=& \frac 1 2 \Big\{  \big(\big[\mathds{1} +L^* L+\a \h L^* \h L \big]U_{h\t} ,U_{h\t}  \big)_{L^2_{{\mathbb F}}(0,T; \dbL^2)}  
     +2 \big(\big[L^* \G+\a \h L^*\h\G \big] {\cR_h X_0} +L^*f+\a \h L^*\h f ,U_{h\t} \big)_{L^2_{{\mathbb F}}(0,T; \dbL^2)}\\
   &\q  +\Big[ \big( \G {\cR_h X_0}+f, \G {\cR_h X_0}+f\big)_{L^2_{{\mathbb F}}(0,T; \dbL^2)}
   + \a \big( \h\G {\cR_h X_0}+\h f, \h\G \cR_h X_{0}+\h f \big)_{L^2_{\mf_T}(\O;\dbL^2)} \Big]     \Big\}\\
=:& \frac 1 2\Big[ \big( NU_{h\t} ,U_{h\t} \big)_{L^2_{{\mathbb F}}(0,T; \dbL^2)} +2 \big( H({\cR_h X_0,f}),U_{h\t}\big)_{L^2_{{\mathbb F}}(0,T; \dbL^2)}
    + M(\cR_hX_{0},f)  \Big]\, .    
\eeq
We use this re-writing of $\cJ_{\t}(X_{h\t} ,U_{h\t} )$  which involves
mappings $N$ and $H$, and where the last term does not depend on $U_{h\t}$,
to now involve the optimality condition for $U^*_{h\t}$: since
$N  =\mathds{1} +L^* L+\a \h L^* \h L $  is positive definite, there exists a unique $U^*_{h\t} \in {\mathbb U}_{h\t}$ such that
\beq
NU^*_{h\t} +H({\cR_h X_0},f ) =0\, .
\eeq
Therefore, for any $U_{h\t} \in {\mathbb U}_{h\t}$ such that $U_{h\t} \neq U^*_{h\t}$, 
\beq
\cJ_{\t}(X_{h\t},U_{h\t})-\cJ_{\t}(X^*_{h\t}, U^*_{h\t})
=\frac 1 2 \big( N(U_{h\t}-U^*_{h\t}), U_{h\t}-U^*_{h\t}\big)_{L^2_{{\mathbb F}}(0,T; \dbL^2)}
>0\, ,
\eeq
which means that $U^*_{h\t}$ is the unique optimal control, and $(X^*_{h\t},U^*_{h\t})$ is the unique optimal pair.

Finally,  by the definition of $N, H,L^*,\h L^*$,  and properties \eqref{w1003e01} and
 \eqref{rep1}, we can get
 \beq
0=N U^*_{h\t} +H({\cR_h X_0})
=U^*_{h\t} +\big[L^* X^*_{h\t} +\a \h L^*  X^*_{h\t}(T) \big]\, . \\
\eeq
By the definition of $K_{h\t}$ in \rf{w1212e3}$_2$, we can arrive at 
\beq
L^* X^*_{h\t} +\a \h L^*  X^*_{h\t}(T)=-K_{h\t}X^*_{h\t}\,.
\eeq
Then \eqref{w1003e12} can be deduced by these two equalities. That completes the proof.
\end{proof}


In Remark \ref{w829r1} we already mentioned that $K_{h\t}X^*_{h\t}$  in 
\eqref{w1212e3}$_2$ neither solves the temporal discretization of  {\bf BSPDE}$_h$ \eqref{bshe1a} by the explicit Euler, nor by the implicit Euler method.
The following two lemmas study the difference between $K_{h\t}X^*_{h\t}$ and the temporal discretization of {\bf BSPDE}$_h$ \eqref{bshe1a} by
the implicit Euler method. 
\bl{w229l1}
Suppose that $(Y_0,\, Z_0)$ solves the following {\bf BSDE}$_h$:
\bel{w229e7}
\lt\{
\bal
&\rd Y_0(t) =  \bigl[-\D_h Y_0\bigl(\nu(t)\bigr)- \bar Z_0 \bigl(\nu(t) \bigr)
+ X^*_{h\t}\bigl(\mu(t)\bigr) \bigr]\rd t +\ Z_0(t)\rd W(t)  \q \forall~ t\in [0,T]\, ,\\
&Y_0(t_N)=-\a X^*_{h\t}(T) \, ,
\eal
\rt.
\ee
where $X^*_{h\t}$ is the optimal state of problem {\bf SLQ$_{h\t}$}, $\mu(\cd),\, \nu(\cd)$ are defined in \rf{w827e1}, 
and  $\bar Z_0$ is a piecewise constant process
 which is defined by
\beq
\bar Z_0(t)=\frac 1 \t \me\Big[\int_{t_n}^{t_{n+1}}  Z_0(s)\rd s\Big| \mf_{t_n}\Big]\qquad \forall\, t\in [t_n,t_{n+1})\,, \qquad  n=0,1,\cds,N-1\, .
\eeq
 Then, %
\beq
\bal
Y_0(t_j)
=&-\a \me \Big[ A_0^{N-j}\prod_{k={j+1}}^N \lt(1+\D_kW\rt)  X^*_{h\t}(T) \Big|\mf_{t_j} \Big]\\
 &     - \t \me \Big[ \sum_{k=j+1}^NA_0^{k-j}\prod_{i=j+1}^{k}\lt(1+\D_iW\rt)  X^*_{h\t}(t_k) \Big|\mf_{t_j} \Big] \qquad j=0,1,\cds,N-1\, .\\
\eal
\eeq
\el
Note that $(Y_0,\bar Z_0)$ is the numerical solution to {\bf BSPDE}$_h$ \eqref{bshe1a} by the implicit Euler method;
 see e.g.~\cite{Bender-Denk07}.
\begin{proof}
By \eqref{w229e7}, for any $n=0,1,\cds,N-1$,
$$
\t \bar Z_0(t_n)
=\me\Big[ \bigl( Y_0(t_{n+1})-\t X^*_{h\t}(t_{n+1}) \bigr) \D_{n+1}W\Big|\mf_{t_n} \Big]\,,
$$
which yields the desired results.
\end{proof}

The solution $(Y_0,\,Z_0)$ in \rf{w229e7} depends on $X^*_{h\t}$. Hence, in what follows, we may write it
in the form $\bigl(Y_0(\cd;X^*_{h\t}),  \,Z_0(\cd;X^*_{h\t})\bigr)$. 

Similar to $\cS,\,\cS_h$, we can define difference equation \rf{w1003e7}'s solution operator by 
\beq
\cS_{h\t}: \dbU_{h\t}\to \dbX_{h\t}\,.
\eeq
In the next lemma, we estimate the 
difference between $(K_{h\t}\cS_{h\t}(\Pi_\t U^*_h))(\cd)$,
which was introduced in \rf{w1212e3}$_2$, and $Y_0\bigl(\cd;\cS_{h\t}(\Pi_\t U^*_h)\bigr)$, which is crucial in proving rates of
convergence for the temporal discretization.

\bl{w301l1} Suppose that $(Y_0,Z_0)$ solves \eqref{w229e7}. Then
it holds that
\beq
\max_{0\leq j\leq N-1}\me\Big[ \Vert Y_0\bigl(t_j; \cS_{h\t}(\Pi_\t U^*_h)\bigr)-\bigl(K_{h\t}\cS_{h\t}(\Pi_\t U^*_h)\bigr)(t_j)
\Vert_{\dbL^2 }^2\Big]\leq C\t \,,
\eeq
where $\Pi_\t$ is defined in \eqref{w828e1}, where $K_{h\t}$ is given by \eqref{w1212e3}$_2$, and $C$ is independent of $h$ and $\t$.
\el

\begin{proof} By triangular inequality, it holds that, for any $j=0,1,\cds,N-1$,
\bel{w301e15}
\bal
&\me\big[ \Vert Y_0(t_j; \cS_{h\t}(\Pi_\t U^*_h))-(K_{h\t}\cS_{h\t}(\Pi_\t U^*_h))(t_j)
\Vert_{\dbL^2}^2\big]\\
&\q\leq 2\Biggl(\a^2\me\Big[\Big\| A_0^{N-j}\Big(\prod_{k={j+1}}^N-\prod_{k={j+2}}^N\Big) \lt(1+\D_kW\rt) \cS_{h\t}(T;\Pi_\t U^*_h) \Big\|_{\dbL^2}^2\Big]\\
&\qq\qq+\me\Big[\Big\| \t\sum_{k=j+1}^NA_0^{k-j}\Big(\prod_{i=j+1}^{k}-\prod_{i=j+2}^{k}\Big)\lt(1+\D_iW\rt)  \cS_{h\t}(t_k;\Pi_\t U^*_h) \Big\|_{\dbL^2}^2\Big]\Biggr)\\
&\q=:2\lt(I_{j1}+I_{j2}\rt)\,.
\eal
\ee

We only present the estimate of $I_{j2}$, since $I_{j1}$ can be bounded in a similar vein. 
Applying \eqref{w1212e3}$_1$, we can get
\bel{w301e16}
\bal
I_{j2}\leq& \t^2 N\sum_{k=j+1}^{N} \me\Big[\Big\| \prod_{i=j+2}^{k}\lt(1+\D_iW\rt)\D_{j+1}W
 \cS_{h\t}(t_k;\Pi_\t U^*_h) \Big\|_{\dbL^2}^2\Big]\\
\leq & 3 \t \sum_{k=j+1}^{N}\Biggl(\me\Big[\Big\| \prod_{i=j+2}^{k}\lt(1+\D_iW\rt)\D_{j+1}W
\prod_{l=1}^k\lt(1+\D_l W\rt)X_{h\t}({0}) \Big\|_{\dbL^2}^2\Big]\\
 &\q+\t^2 N\sum_{l=0}^{k-1}\me\Big[\Big\| \prod_{i=j+2}^{k}\lt(1+\D_iW\rt)\D_{j+1}W
\prod_{m=l+2}^k\lt(1+\D_mW\rt) U^*_{h}(t_l) \Big\|_{\dbL^2}^2\Big]\\
&\q + \sum_{l=0}^{k-1}\me\Big[\Big\| \prod_{i=j+2}^{k}\lt(1+\D_iW\rt)\D_{j+1}W
\prod_{m=l+2}^k\lt(1+\D_mW\rt) \cR_{h}\si(t_l)\D_{l+1}W \Big\|_{\dbL^2}^2 \Big]\Biggr) \,.\\
\eal
\ee
For the last term in \rf{w301e16}, we use the fact that  $\{\D_jW\}_{j=1}^N$ is a sequence of mutually independent random variables, 
such that for $0\leq l_1,l_2\leq k-1$, and $l_1\neq l_2$,
\beq
\me\Big[\Big(\prod_{m=l_1+2}^k\lt(1+\D_mW\rt)\D_{l_1+1}W\Big)
\times \Big(\prod_{m=l_2+2}^k\lt(1+\D_mW\rt)\D_{l_2+1}W\Big)\Big]=0\,.
\eeq
Also by the mutual independence of $\{\D_jW\}_{j=1}^N$, and since $X_{h\t}(0)$  
is deterministic, we can arrive at
\bel{w302e1}
\bal
&\me\Big[\Big\| \prod_{i=j+2}^{k}\lt(1+\D_iW\rt)\D_{j+1}W
\prod_{i=1}^k\lt(1+\D_i W\rt)X_{h\t}({0}) \Big\|_{\dbL^2}^2\Big]\\
&\q=\prod_{i=j+2}^k\me\big[(1+\D_iW)^4\big] \me\big[\big(\D_{j+1}W+\D_{j+1}^2W\big)^2\big]
  \prod_{i=1}^j\me\big[(1+\D_iW)^2\big]\lt\| X_{h\t}(0)\rt\|_{\dbL^2}^2\\
&\q=  \lt(\t+3\t^2\rt)\prod_{i=j+2}^k (1+4\t+3 \t^2)
  \prod_{i=1}^j (1+\t)\lt\| X_{h\t}(0)\rt\|_{\dbL^2}^2\\
&\q \leq C \t  \lt\| X(0)\rt\|_{\dbH_0^1}^2\,,
\eal
\ee
where $\D_{j+1}^2W:=\lt(\D_{j+1}W \rt)^2$.

In the sequel,  we tend to estimate $ \me\Big[\Big\| \prod_{i=j+2}^{k}\lt(1+\D_iW\rt)\D_{j+1}W
\prod_{m=l+2}^k\lt(1+\D_mW\rt) U^*_{h}(t_l) \Big\|_{\dbL^2}^2\Big]$ and $\me\Big[\Big\| \prod_{i=j+2}^{k}\lt(1+\D_iW\rt)\D_{j+1}W
\prod_{m=l+2}^k\lt(1+\D_mW\rt) \cR_{h}\si(t_l)\D_{l+1}W \Big\|_{\dbL^2}^2\Big]$
under the 
following two cases.

\no{\bf Case I.}  $l\leq j$. In this case, using the same trick as that in \eqref{w302e1}, 
we can see that
\bel{w302e2}
\bal
&\me\Big[\Big\| \prod_{i=j+2}^{k}\lt(1+\D_iW\rt)\D_{j+1}W
\prod_{m=l+2}^k\lt(1+\D_mW\rt) U^*_{h}(t_l) \Big\|_{\dbL^2}^2\Big]\\
&\q \leq \prod_{i=j+2}^k\me\big[(1+\D_iW)^4\big]\me\big[(\D_{j+1}W+\D_{j+1}^2W )^2\big]
  \prod_{m=l+2}^j\me\big[(1+\D_mW)^2\big]\me\big[\|U^*_h(t_l)\|_{\dbL^2}^2\big]\\
&\q\leq C\t \sup_{t\in[0,T]}\me\big[\|U^*_h(t)\|_{\dbL^2}^2\big]\,,
\eal
\ee
and
\bel{w1016e2}
\bal
&\me\Big[\Big\| \prod_{i=j+2}^{k}\lt(1+\D_iW\rt)\D_{j+1}W
\prod_{m=l+2}^k\lt(1+\D_mW\rt) \cR_{h}\si(t_l)\D_{l+1}W \Big\|_{\dbL^2}^2\Big]\\
&\q\leq
\lt\{
\bal
&\prod_{i=j+2}^k\me\big[(1+\D_iW)^4\big]\me\big[(\D_{j+1}W+\D_{j+1}^2W )^2\big]\\
 &\qq\qq\qq\times \prod_{m=l+2}^j\me\big[(1+\D_mW)^2\big]\me\big[\D_{l+1}^2W\big]\|\si(t_l)\|_{\dbH_0^1}^2 \q  l<j,\\
&\prod_{i=j+2}^k\me\big[(1+\D_iW)^4\big]\me\big[(\D_{j+1}W )^4\big]
\|\si(t_l)\|_{\dbH_0^1}^2 \q  l=j,\\
\eal
\rt.  \\
&\q\leq C\t^2 \|\si\|_{C([0,T];\dbH_0^1)}^2\,.
\eal
\ee

\no{\bf Case II.}  $l> j$. Still applying the mutual independence of $\{\D_jW\}_{j=1}^N$, we can deduce that
\bel{w302e3}
\bal
&\me\Big[\Big\| \prod_{i=j+2}^{k}\lt(1+\D_iW\rt)\D_{j+1}W
\prod_{m=l+2}^k\lt(1+\D_mW\rt) U^*_{h}(t_l) \Big\|_{\dbL^2}^2\Big]\\
&\q\leq \prod_{i=l+2}^k\me\big[(1+\D_i W)^4\big]\me\big[(1+\D_{l+1}W )^2\big]\\
&\qq\qq\times \Big\{ \prod_{m=j+2}^l\me\big[(1+\D_m W)^4\big]\me\big[(\D_{j+1}W)^4\big]\me\big[\|U^*_h(t_l)\|_{\dbL^2}^4\big]\Big\}^{1/2}\\
&\q \leq C\t  \Big(\sup_{t\in[0,T]} \me\big[\|U^*_h(t)\|_{\dbL^2}^4\big] \Big)^{1/2}\, , 
\eal
\ee
and
\bel{w1016e3}
\bal
&\me\Big[\Big\| \prod_{i=j+2}^{k}\lt(1+\D_iW\rt)\D_{j+1}W
\prod_{m=l+2}^k\lt(1+\D_mW\rt)\Pi_{h}\si(t_l)\D_{l+1}W \Big\|_{\dbL^2}^2\Big]\\
&\q=\prod_{i=l+2}^k\me\big[(1+\D_iW)^4\big] \me\big[(\D_{l+1}W+\D_{l+1}^2W )^2\big]\\
&\qq\qq\times  \prod_{m=j+2}^l\me\big[(1+\D_mW)^2\big]\me\big[\D_{j+1}^2W\big]\|\cR_h\si(t_l)\|_{\dbL^2}^2 \\
&\q \leq C\t^2 \|\si\|_{C([0,T];\dbH_0^1)}^2\, . 
\eal
\ee

In the above both cases, utilizing \eqref{w826e1}$_6$ of Lemma \ref{w229l2},  and combining with 
\rf{w301e16}, \rf{w302e1}, we can bound $I_{j2}$ by $C\t$, 
and then prove the desired
result.
\end{proof}

\bl{w308l1}
Suppose that 
$(Y_h,Z_h)$, $(Y_0,Z_0)$ solve \eqref{bshe1a}, \eqref{w229e7} respectively. Then it holds that
\beq
 \max_{0\leq n\leq N}\me \big[\|Y_h(t_n;\cS_{h\t}(\Pi_\t U^*_h))-Y_0(t_n; \cS_{h\t}(\Pi_\t U^*_h))\|_{\dbL^2}^2\big] 
\leq   C\t\, .
\eeq
\el

\begin{proof}
We use a standard argument for the implicit method to solve BSDEs (see e.g.~\cite{ZhangJF04, Wang16}); since we need to deal with a family of BSDEs, 
which is obtained via a finite element discretization of  {\bf BSPDE} and thus depends on the parameter $h$, for the completeness, we provide the proof. 

For simplicity, we write $$\Bigl(Y_h\bigl(\cd;\cS_{h\t}(\Pi_\t U_h^*)\bigr), Z_h\bigl(\cd;\cS_{h\t}(\Pi_\t U_h^*)\bigr)\Bigr)   \quad \mbox{resp.} \quad \Bigl(Y_0\bigl(\cd;\cS_{h\t}(\Pi_\t U_h^*)\bigr), Z_0\bigl(\cd;\cS_{h\t}(\Pi_\t U_h^*)\bigr)\Bigr)$$ as $\bigl(\wt Y_h (\cd),\wt Z_h(\cd)\bigr)$ resp.~$\bigl(\wt Y_0(\cd),\wt Z_0(\cd) \bigr)$
in the proof.
Set $ e_Y^n=\wt Y_h(t_n)-\wt Y_0(t_n)$. Subtracting \eqref{w229e7} from 
\eqref{bshe1a} by changing $X^*_h$ to $\cS_{h\t}(\Pi_\t U^*_h)$, we can get
\bel{w1124e1}
\bal
&(\mathds{1}-\t\D_h) e_Y^n +\int_{t_n}^{t_{n+1}} \bigl(\wt Z_h(s)- \wt Z_0(s)\bigr)\rd W(s)\\
&\q= e_Y^{n+1} +\int_{t_n}^{t_{n+1}}\Big[ \D_h \bigl(\wt Y_h(t_n)-\wt Y_h(s) \bigr)
   +\frac 1 \t \me\big[\int_{t_n}^{t_{n+1}} \wt Z_0(t)\rd t\big|\mf_{t_n} \big]- \wt Z_h(s)\\
  &\qq\qq\qq  \qquad +\Bigl(\cS_{h\t}(t_n, \Pi_\t U^*_h)-\cS_{h\t}(t_{n+1},\Pi_\t U^*_h) \Bigr) \Big]\rd s\, .
\eal
\ee
Following the same procedure as in the proof of \cite[Theorem 4.1]{Wang16}, we get for all $\varepsilon >0$,
\beq
\bal
& \me\big[\|(\mathds{1}-\t\D_h ) e_Y^n  \|^2_{\dbL^2}\big]+\me\Big[\int_{t_n}^{t_{n+1}}\big\| \wt Z_h(s)-\wt Z_0(s)\big\|_{\dbL^2}^2\rd s\Big]\\
&\q \leq  (1+5\e)\me\big[\| e_Y^{n+1}  \|^2_{\dbL^2}\big]
  +(5+1/\e)\Biggl( \t\me\Big[\int_{t_n}^{t_{n+1}} \big\| \D_h(\wt Y_h(s)-\wt Y_h(t_n))\big\|_{\dbL^2}^2 \rd s\Big]\\
 &\qq+2\t \me\Big[ \int_{t_n}^{t_{n+1}} \big\| \wt Z_h(s)-\wt Z_h(t_n)\big\|_{\dbL^2}^2\rd s \Big]
   +\t \me\Big[\int_{t_n}^{t_{n+1}} \big\| \wt Z_h(s)- \wt Z_0(s)\big\|_{\dbL^2}^2\rd s\Big]\\
 &\qq+\t^2 C \me \big[\|\cS_{h\t}(t_n, \Pi_\t U^*_h)-\cS_{h\t}(t_{n+1},\Pi_\t U^*_h) \|_{\dbL^2}^2\big]   \Biggr)\,.
\eal
\eeq
By taking $\e\geq \frac{\t}{1-5\t}$ (for example we may choose $\t$ small enough such that
$ 5\t\leq  1/2 $ and take $\ds \e=2\t$), we will see that
\beq
\bal
\me\big[\| e_Y^n  \|^2_{\dbL^2}\big]
\leq &(1+10\t) \me\big[\| e_Y^{n+1}  \|^2_{\dbL^2}\big]\\
 &+\Big(5+\frac 1{2\t}\Big)\Biggl( \t\me\Big[\int_{t_n}^{t_{n+1}} \big\| \D_h\bigl(\wt Y_h(s)-\wt Y_h(t_n)\bigr)\big\|_{\dbL^2}^2
      +2 \big\| \wt Z_h(s)-\wt Z_h(t_n)\big\|_{\dbL^2}^2 \rd s\Big]\\
 &
   \qq\qq\qq+\t^2 C \me \big[\|\cS_{h\t}(t_n, \Pi_\t U^*_h)-\cS_{h\t}(t_{n+1},\Pi_\t U^*_h) \|_{\dbL^2}^2 \big]  \Biggr)\,.
\eal
\eeq
Then the discrete Gronwall's inequality leads to
\bel{w308e2}
\bal
\max_{0\leq n\leq N-1}\me\big[\| e_Y^n  \|^2_{\dbL^2}\big]
&\leq  C   \sum_{n=0}^{N-1}\Big\{\me\Big[\int_{t_n}^{t_{n+1}} \big\| \D_h\bigl(\wt Y_h(s)-\wt Y_h(t_n)\bigr)\big\|_{\dbL^2}^2
  +\big\| \wt Z_h(s)-\wt Z_h(t_n)\big\|_{\dbL^2}^2 \rd s \Big]\\
 &\qq\qq+\t \me \big[\|\cS_{h\t}(t_n, \Pi_\t U^*_h)-\cS_{h\t}(t_{n+1},\Pi_\t U^*_h) \|_{\dbL^2}^2\big]   \Big\}\\
 &=:C\lt(I_1+I_2+I_3 \rt)  \,.
\eal
\ee
In what follows, we estimate the three terms on the right side of \eqref{w308e2}. 

For $I_1$,   by {\bf BSPDE}$_h$ \eqref{bshe1a}  (changing $X^*_h$ to $\cS_{h\t}(\Pi_\t U^*_h)$) and It\^o's formula to
$\|\nb\D \wt Y\|_{\dbL^2}^2$, we estimate
\bel{w906e1}
\bal
 I_1\leq& C\t \me\Big[\int_0^T  \|\D_h^2 \wt Y_h(t) \|_{\dbL^2}^2
   +\|\D_h \wt Z_h(t)\|_{\dbL^2}^2 \rd t \Big]
   +C\t  \sup_{0\leq n\leq N }\me[ \|\D_h \cS_{h\t}(t_n;\Pi_\t U^*_h)\|_{\dbL^2}^2] \\
\leq & C\t 
  \Big(\me\big[ \|\nb \D_h  \cS_{h\t}(T;\Pi_\t U^*_h)\|_{\dbL^2}^2\big]
   + \sup_{0\leq n\leq N}\me \big[\|\D_h \cS_{h\t}(t_n;\Pi_\t U^*_h)\|_{\dbL^2}^2\big] \Big)\,.\\
\eal
\ee
In the below, we tend to estimate two terms on the right side of \rf{w906e1}.
By Lemma \ref{w1011l1} and the following calculation:
%
\bel{w309e2}
\bal
&\me\big[\|\D_h\cS_{h\t}(t_{n};\Pi_\t U^*_h) \|_{\dbL^2}^2\big]\\
\leq & 3 \prod_{j=1}^n\lt(1+\t\rt)\lt\|\D_h X_{h\t}({0}) \rt\|_{\dbL^2}^2
 +3\t^2 N \sum_{j=0}^{n-1} \prod_{k=j+2}^n\lt(1+\t\rt)\me\big[\|\D_hU_{h}^*(t_j) \|_{\dbL^2}^2\big]\\
& +3\t  \sum_{j=0}^{n-1} \prod_{k=j+2}^n(1+\t)\|\D_h \cR_h\si (t_j) \|_{\dbL^2}^2\\
\leq& C  \lt[\|X_0\|_{\dbH_2}^2+\|\si\|^2_{C([0,T];\dbH^2)}\rt]\, ,
\eal
\ee
we can see that
 the second term on the right side of \rf{w906e1} is bounded. In the below,
we estimate the first term.
By \eqref{w1212e3}$_1$ (with $U_{h\t}=\Pi_\t U^*_h$), for any $k=0,1,\cds,N$, 
\bel{w922e2}
\bal
\me\Big[\Big\|\nb\D_hA_0^N\prod_{j=1}^N(1+\D_jW)X_{h\t}(0)\Big\|_{\dbL^2}^2\Big]
=(1+\t)^N \|\nb\D_h \h X_{h\t}(T)\|_{\dbL^2}^2\,,
\eal
\ee
where $\h X_{h\t}(t_k)=A_0^k X_{h\t}(0)$, for $k=0,1,\cds,N$ solves
\bel{w908e5}
\lt\{
\bal
&\h X_{h\t}(t_{k+1})-\h X_{h\t}(t_k)=\t \D_h \h X_{h\t}(t_{k+1})\q k=0,1,\cds,N-1\,,\\
&\h X_{h\t}(0)=X_{h\t}(0)\,.
\eal
\rt.
\ee
Multiplying \rf{w908e5}$_1$ by $\D_h^3\h X_{h\t}(t_{k+1})$ then leads to 
\bel{w922e3}
\bal
&\|\nb\D_h \h X_{h\t}(t_{k+1})\|_{\dbL^2}^2=\|\nb\D_h A_0\h X_{h\t}(t_{k})\|_{\dbL^2}^2\\
=&\|\nb\D_h \h X_{h\t}(t_k)\|_{\dbL^2}^2
   -\|\nb\D_h \big(\h X_{h\t}(t_{k+1})-\h X_{h\t}(t_k)\big)\|_{\dbL^2}^2
   -\t \|\D_h^2 \h X_{h\t}(t_{k+1})\|_{\dbL^2}^2\\
\leq &    \|\nb\D_h \h X_{h\t}(t_k)\|_{\dbL^2}^2\,,
\eal
\ee
which, together with \rf{w922e2}, leads to
\bel{w922e4}
\bal
\me\Big[\Big\|\nb\D_hA_0^N\prod_{j=1}^N(1+\D_jW)X_{h\t}(0)\Big\|_{\dbL^2}^2\Big]
\leq C \|\nb\D_h  X_{h\t}(0)\|_{\dbL^2}^2\,.
\eal
\ee
With a similar trick, we can show that
\bel{w1016e4}
\bal
\me\Big[\Big\|\nb\D_hA_0^{n-j}\prod_{k=j+2}^n(1+\D_kW)\cR_h\si(t_j)\D_{j+1}W\Big\|_{\dbL^2}^2\Big]
\leq C \|\nb\D_h \cR_h\si(t_j)\|_{\dbL^2}^2
\leq C\|\si(t_j)\|_{\dbH^3}^2\,.
\eal
\ee

On the other side, by sapplying the maximum condition \rf{pontr1a}, estimate \rf{w922e3}, Lemma \ref{w1011l1} and  It\^o's formula to $\|\nb\D_hA_0^kY_h\|_{\dbL^2}^2$, we can find that
\bel{w922e1}
\bal
&\sup_{t\in [0,T]}\me\big[\|\nb\D_hA_0^k U^*_h(t)\|_{\dbL^2}^2\big]
+\me\Big[\int_0^T\|\D_h^2A_0^k U^*_h(t)\|_{\dbL^2}^2+\|\nb\D_hA_0^k Z_h(t)\|_{\dbL^2}^2\rd t \Big]\\
&\leq C\Big[\me\big[\|\nb\D_hA_0^kX^*_h(T)\|_{\dbL^2}^2\big]+\|\nb\D_h X_h^*(0)\|_{\dbL^2}^2\Big]\\
&\leq C\Big[\me\big[\|\nb\D_h X^*_h(T)\|_{\dbL^2}^2\big]+\|\nb \D_hX_h^*(0)\|_{\dbL^2}^2\Big]\\
&\leq C\lt[ \|X_0\|_{\dbH^3}^2+\|\si\|^2_{L^2(0,T;\dbH^3)}\rt]\,.
\eal
\ee

Combining with \rf{w922e4}, \rf{w1016e4} and \rf{w922e1}, we can deduce that
\beq
\me \big[\|\nb \D_h  \cS_{h\t}(T;\Pi_\t U^*_h)\|_{\dbL^2}^2\big]
\leq C\lt[ \|X_0\|_{\dbH^3}^2+\|\si\|^2_{C([0,T];\dbH^3)}\rt]\,.
\eeq
Therefore, by \eqref{w906e1}, we arrive at
\beq
I_1\leq C\lt[ \|X_0\|_{\dbH^3}^2+\|\si\|^2_{C([0,T];\dbH^3)}\rt]\t \,.
\eeq

For $I_2$, we still can prove that  
\bel{w309e4}
I_2\leq C\Big[\|X_0\|_{\dbH_0^1}^2+\|\si\|_{C([0,T];\dbH_0^1)}^2\Big]\t \,.
\ee
We postpone the proof to  Lemma \ref{w1024l1}.

For $I_3$,
\eqref{w1003e7} (with $U_{h\t}=\Pi_\t U^*_h$) yields to
\bel{w308e3}
\bal
&\me \big[\|\cS_{h\t}(t_n; \Pi_\t U^*_h)-\cS_{h\t}(t_{n+1}; \Pi_\t U^*_h) \|_{\dbL^2}^2\big]\\
&\q=\me \Big[\| \tau \big[\D_h \cS_{h\t}(t_{n+1}; \Pi_\t U^*_h) +
{ U_{h}^*(t_n)} \big] +\lt[ \cS_{h\t}(t_{n}; \Pi_\t U^*_h) +\cR_h\si(t_n)\rt] \D_{n+1}W  \|_{\dbL^2}^2\Big] \\
&\q\leq  C\t \Big[\max_{0\leq n\leq N}\me\big[\|\cS_{h\t}(t_{n}; \Pi_\t U^*_h) \|_{\dbL^2}^2\big]
+ \max_{0\leq n\leq N}\me\big[\|\D_h\cS_{h\t}(t_{n}; \Pi_\t U^*_h) \|_{\dbL^2}^2\big]\\
&\qq\q +\sup_{t\in[0,T]}\me\big[\| U^*_h(t) \|_{\dbL^2}^2\big]+\sup_{t\in[0,T]}\lt\| \si(t) \rt\|_{\dbH_0^1}^2\Big]\\
&\q\leq  C\lt[\|X_0\|_{\dbH_2}^2+\|\si\|^2_{C([0,T];\dbH^2)}\rt]\t\,.
\eal
\ee
%
\eqref{w309e2}, together with \eqref{w826e1}$_5$ and \eqref{w308e3} then leads to
\begin{eqnarray}\nonumber
I_3 &\leq& C \max_{0\leq n\leq N-1}
\me \big[\|\cS_{h\t}(t_n; \Pi_\t U^*_h)-\cS_{h\t}(t_{n+1}; \Pi_\t U^*_h) \|_{\dbL^2}^2\big] \\ \label{308e6}
&\leq&  C \lt[\|X_0\|_{\dbH_2}^2+\|\si\|^2_{C([0,T];\dbH^2)}\rt] \t \, .
\end{eqnarray}
That completes the proof.
\end{proof}

\bl{w1024l1}
Suppose that 
$(Y_h,Z_h)$ solves 
\eqref{bshe1a}. Then  for all $t\in [0,T]$,
\bel{w1024e9}
\bal
&\me\big[\|Z_h(t; \cS_{h\t}(\Pi_\t U^*_h))-Z_h(\nu(t); \cS_{h\t}(\Pi_\t U^*_h))\|_{\dbL^2}^2 \big]
\leq   C|t-\nu(t)|  \Big(\|X_0\|_{\dbH_0^1}^2+\|\si\|_{C([0,T];\dbH_0^1)}^2\Big)  \, ,
\eal
\ee
where $ \nu(\cd)$ is defined in \rf{w827e1}
\el

\begin{proof} The proof is long. Hence we divide it into three steps.

\no{\bf 1)} 
We claim that, for any $n=0,1,\cds,N$, $\cS_{h\t}(t_n;\Pi_\t U^*_h)\in \dbD^{2,2}(\dbL^2)$.

Indeed, by  \rf{w1024e2}, we know that 
\beq
\cS_{h\t}(t_n;\Pi_\t U^*_h)
=&A_0^n\prod_{j=1}^n\lt(1+\D_j W\rt)X_{h\t}({0})
+\t \sum_{j=0}^{n-1}A_0^{n-j}\prod_{k=j+2}^n\lt(1+\D_kW\rt) U_{h}^*(t_j)\\
&+\sum_{j=0}^{n-1} A_0^{n-j}\prod_{k=j+2}^{n}\lt(1+\D_kW \rt)\cR_h\si(t_j)\D_{j+1}W\,.
\eeq
In the following, we only prove that the second term on the right-hand side of the above representation is in $\dbD^{2,2}(\dbL^2)$.
The other two terms can be proved in a similar vein.

By maximum condition \rf{pontr1a} and Lemma \ref{w229l3}, 
we know that for any $j=0,1,\cds,N-1$, $U^*_h(t_j)\in \dbD^{2,2}(\dbL^2)$. By the chain rule, for any $\th,\mu \in [0,T]$
without loss of generality, suppose that $\th\in [t_l,t_{l+1}),\, \mu\in [t_m,t_{m+1})$. Then, by the fact that 
$D_\th (1+\D_kW)=\d_{lk},\,D_\mu (1+\D_kW)=\d_{mk}$, where $\d_{lk},\,\d_{mk}$ are Kronecker delta functions, we have
\beq
&D_\th \Big(\t \sum_{j=0}^{n-1}A_0^{n-j}\prod_{k=j+2}^n\lt(1+\D_kW\rt) U_{h}^*(t_j)\Big)\\
&=\t \sum_{j=0}^{n-1}A_0^{n-j}\prod_{\scriptstyle k=j+2 \atop \scriptstyle k\neq l}^n \lt(1+\D_kW\rt) U_{h}^*(t_j)
  +\t \sum_{j=0}^{n-1}A_0^{n-j}\prod_{k=j+2}^n\lt(1+\D_kW\rt) D_\th U_{h}^*(t_j)\,,
\eeq
\beq
&D_\mu D_\th \Big(\t \sum_{j=0}^{n-1}A_0^{n-j}\prod_{k=j+2}^n\lt(1+\D_kW\rt) U_{h}^*(t_j)\Big)\\
&=\t \sum_{j=0}^{n-1}A_0^{n-j}\d_{ml}\prod_{\scriptstyle k=j+2 \atop \scriptstyle k\neq l,k\neq m}^n \lt(1+\D_kW\rt) U_{h}^*(t_j)
 +\t \sum_{j=0}^{n-1}A_0^{n-j}\prod_{\scriptstyle k=j+2 \atop \scriptstyle k\neq l}^n \lt(1+\D_kW\rt) D_\mu U_{h}^*(t_j)\\
&\q  +\t \sum_{j=0}^{n-1}A_0^{n-j}\prod_{\scriptstyle k=j+2 \atop \scriptstyle k\neq m}^n\lt(1+\D_kW\rt) D_\th U_{h}^*(t_j)
 +\t \sum_{j=0}^{n-1}A_0^{n-j}\prod_{k=j+2}^n\lt(1+\D_kW\rt) D_\mu D_\th U_{h}^*(t_j) \,,
\eeq
and subsequently, by \rf{w301e1}$_1$ and \rf{w301e1}$_4$ in Lemma  \ref{w229l3}, 
\beq
\bal
&\me\Big[\Big\|D_\mu D_\th \Big(\t \sum_{j=0}^{n-1}A_0^{n-j}\prod_{k=j+2}^n\lt(1+\D_kW\rt) U_{h}^*(t_j)\Big)\Big\|_{\dbL^2}^2\Big]\\
&\leq C\Big[\sup_{0\leq j\leq N-1}\me\big[\|U^*_h(t_j)\|_{\dbL^2}^2\big]+\sup_{0\leq j\leq N-1}\me\big[\|D_\th U^*_h(t_j)\|_{\dbL^2}^2\big]
  +\sup_{0\leq j\leq N-1}\me\big[\|D_\mu D_\th U^*_h(t_j)\|_{\dbL^2}^2\big] \Big]\\
&\leq C \Big[\|X_0\|_{\dbH_0^1}^2+\|\si\|_{C([0,T];\dbL^2)}^2+\|\si\|_{L^2(0,T;\dbH_0^1)}^2\Big]\,,
\eal
\eeq
which leads to $\t \sum_{j=0}^{n-1}A_0^{n-j}\prod_{k=j+2}^n\lt(1+\D_kW\rt) U_{h}^*(t_j)\in \dbD^{2,2}(\dbL^2)$.
Hence, for any $t\in[0,T]$, $\cS_{h\t}(\nu(t);\Pi_\t U^*_h)\in \dbD^{2,2}(\dbL^2)$, and
\bel{w1024e8}
\bal
 \sup_{\mu,\th\in[0,T]}\sup_{\nu(t)\in[\mu\vee \th,T]} \me\big[\|D_\mu D_\th \cS_{h\t}(\nu(t);\Pi_\t U^*_h)\|_{\dbL^2}^2\big]\leq C \Big[\|X_0\|_{\dbH_0^1}^2+\|\si\|_{C([0,T];\dbL^2)}^2+\|\si\|_{L^2(0,T;\dbH_0^1)}^2\Big]\,.
\eal
\ee

\ms
\no{\bf 2)} Applying the same trick to estimate $\me \big[\|\nb \D_h  \cS_{h\t}(T;\Pi_\t U^*_h)\|_{\dbL^2}^2\big]$ 
in the proof of Lemma \ref{w308l1}, we can deduce that for any $\th\in [0,T]$,
\bel{w1024e6}
\bal
\me\big[ \|\nb D_\th \cS_{h\t}(T;\Pi_\t U^*_h)\|_{\dbL^2}^2\big]
\leq C \Big[\|X_0\|_{\dbH_0^1}^2+\|\si\|_{C([0,T];\dbH_0^1)}^2\Big]\,.
\eal
\ee
Also, by the same procedure as that in the proof of \rf{w301e1}$_4$, \rf{w301e1}$_5$ in Lemma \ref{w229l3}, 
thanks to \rf{w1024e8}, we can obtain
\bel{w1024e11}
\bal
&\sup_{\th\in[0,T]}\me\Big[\int_{\th}^T\| D_\th Z_h(t; \cS_{h\t}(\Pi_\t U^*_h))\|_{\dbL^2}^2\rd t\Big]\\
&\q\leq C \sup_{\th\in[0,T]}\sup_{\nu(t)\in[\th,T]}\me\big[\|D_\th \cS_{h\t}(\nu(t);\Pi_\t U^*_h)\|_{\dbL^2}^2\big]\\
&\q\leq C \Big[\|X_0\|_{\dbH_0^1}^2+\|\si\|_{C([0,T];\dbL^2)}^2+\|\si\|_{L^2(0,T;\dbH_0^1)}^2\Big]\,,
\eal
\ee
\bel{w1024e12}
\bal
&\sup_{\mu,\th\in[0,T]}\sup_{t\in[\mu\vee \th,T]}\me\big[\|D_\mu D_\th Y_h(t; \cS_{h\t}(\Pi_\t U^*_h))\|_{\dbL^2}^2\big]\\
&\q\leq C \sup_{\mu,\th\in[0,T]}\sup_{\nu(t)\in[\mu\vee \th,T]}\me\big[\|D_\mu D_\th \cS_{h\t}(\nu(t);\Pi_\t U^*_h)\|_{\dbL^2}^2\big]\\
&\q\leq C \Big[\|X_0\|_{\dbH_0^1}^2+\|\si\|_{C([0,T];\dbL^2)}^2+\|\si\|_{L^2(0,T;\dbH_0^1)}^2\Big]\,,
\eal
\ee
and
\bel{w1024e13}
\bal
&\sup_{\th\in[0,T]}\me\Big[\int_{\th}^T\| \D_h D_\th Y_h(t; \cS_{h\t}(\Pi_\t U^*_h))\|_{\dbL^2}^2\rd t\Big]\\
&\leq C \sup_{\th\in[0,T]}\sup_{\nu(t)\in[\th,T]}\me\big[\|\nb D_\th \cS_{h\t}(\nu(t);\Pi_\t U^*_h)\|_{\dbL^2}^2\big]\\
&\leq C \Big[\|X_0\|_{\dbH_0^1}^2+\|\si\|_{C([0,T];\dbH_0^1)}^2\Big]\,.
\eal
\ee

\ms
\no{\bf 3)} 
Applying the fact that $Z_h(\cd;\cS_{h\t}(\Pi_\t U^*_h))=D_\cd Y_h(\cd;\cS_{h\t}(\Pi_\t U^*_h)) \, \ae$, for any $t\in [t_n,t_{n+1})$, $n=0,1,\cds,N-1$, we arrive at
\bel{w1024e3}
\bal
&\me \bigl[\|Z_h(t;\cS_{h\t}(\Pi_\t U^*_h))-Z_h(\nu(t);\cS_{h\t}(\Pi_\t U^*_h))\|_{\dbL^2}^2\bigr] \\
&\q=\me\bigl[\|D_tY_h(t;\cS_{h\t}(\Pi_\t U^*_h))-D_{t_n}Y_h(t_n;\cS_{h\t}(\Pi_\t U^*_h))\|_{\dbL^2}^2\bigr]\\
&\q \leq 2\me\bigl[\|(D_t-D_{t_n})Y_h(t;\cS_{h\t}(\Pi_\t U^*_h))\|_{\dbL^2}^2\\
&\qq\q +\|D_{t_n}(Y_h(t;\cS_{h\t}(\Pi_\t U^*_h))-Y_h(t_n;\cS_{h\t}(\Pi_\t U^*_h)))\|_{\dbL^2}^2\bigr]\,.
\eal
\ee
For the first term on the right side of \rf{w1024e3}, the fact $t, t_n \in [t_n,t_{n+1})$ leads to 
\beq
(D_t-D_{t_n})\cS_{h\t}(\cd; \Pi_\t U^*_h)\equiv 0 \q \ae\,,
\eeq
and then
\bel{w1024e4}
(D_t-D_{t_n})Y_h(\cd;\cS_{h\t}(\Pi_\t U^*_h))\equiv 0\q \ae\,.
\ee
Similar to \rf{w301e2},  by virtue of \eqref{w1024e8}, \eqref{w1024e11}--\eqref{w1024e8}, we can get
\bel{w1024e5}
\bal
&\me\bigl[\|D_{t_n}(Y_h(t;\cS_{h\t}(\Pi_\t U^*_h))-Y_h(t_n;\cS_{h\t}(\Pi_\t U^*_h)))\|_{\dbL^2}^2\bigr]\\
&\q\leq C\Bigl( \me\Bigl[\int_{t_n}^t\|D_{t_n}D_\th Y_h(\th;\cS_{h\t}(\Pi_\t U^*_h))\|_{\dbL^2}^2\rd \th\Bigr] \\
&\qq\q+(t-{t_n})\me\Bigl[ \int_{t_n}^T \|\D_h D_{t_n}Y_h(\th;\cS_{h\t}(\Pi_\t U^*_h))\|^2_{\dbL^2}+\| D_{t_n}Z_h(\th;\cS_{h\t}(\Pi_\t U^*_h))\|^2_{\dbL^2}\\
&\qq\qq\qq+\|D_{t_n}\cS_{h\t}(\nu(\th);\Pi_\t U^*_h)\|^2_{\dbL^2} \rd \th \Bigr]\Bigr)\\
&\q \leq  C|t-{t_n}|  \Big[\|X_0\|_{\dbH_0^1}^2+\|\si\|_{C([0,T];\dbH_0^1)}^2\Big]\,.
\eal
\ee
Now desired result \rf{w1024e9} can be derived by \rf{w1024e3}--\rf{w1024e5}.
\end{proof}

{
\br{w1124r1}
 For {\bf SPDE} \rf{w1013e1} driven by {\em additive} noise ({\em i.e.}, $\si(t)\rd W(t)$), $Z_h$ does not appear in the drift of
 {\bf BSPDE}$_h$ \rf{bshe1a}. For this case, in order to prove Lemma \ref{w308l1}, we simply multiply by $e_Y^n$ both sides of \rf{w1124e1} and then
 take expectations to settle Lemma \ref{w308l1}. For our problem, however, this
 approach fails due to the presence of $Z_h$ in the drift.
 
\er
}

We are now ready to verify rates of convergence for the solution to problem {\bf SLQ$_{h\t}$}; it is as in Section \ref{rate-1} that the reduced cost functional
$\h\cJ_{h\t}: {\mathbb U}_{h\t} \rightarrow {\mathbb R}$ is used, which is defined via
\beq
\h\cJ_{h\t}(U_{h\t})=\cJ_{\t} \bigl(\cS_{h\t}(U_{h\t}),U_{h\t}\bigr)\, ,
\eeq
where $\cS_{h\t}: {\mathbb U}_{h\tau} \rightarrow {\mathbb X}_{h\tau}$ is the solution operator to the forward equation \eqref{w1212e3}$_1$.

\begin{proof} [\bf {Proof of Theorem \ref{rate2}}]

We divide the proof into two steps.

\no {\bf 1)}
We  follow the argumentation in the proof of Theorem \ref{rate1}.
For every $U_{h\t},\,R_{h\t} \in {\mathbb U}_{h\t}$, the first Fr\'echet derivative $D\h\cJ_{h\t}(U_{h\t})$, and the 
second Fr\'echet derivative
$D^2\h\cJ_{h\t}(U_{h\t})$ satisfy
\bel{w1206e3}
\bal
& D\h\cJ_{h\t}\lt(U_{h\t}\rt)=U_{h\t}-K_{h\t}\cS_{h\t}\lt(U_{h\t}\rt)  \, ,\\
&  \big( D^2 \widehat{\mathcal J}_{h\t}(U_{h\t}) R_{h\t},R_{h\t} \big)_{L^2(0,T; \dbL^2)} 
   \geq \Vert R_{h\t}\Vert^2_{L^2(0,T; \dbL^2)}\, .
\eal
\ee
By putting $R_{h\t} = U^*_{h\t} - \Pi_\t U^*_h$ in \eqref{w1206e3}, and applying the fact
$D\widehat{\mathcal J}_{h\t}(U^*_{h\t})=D\h\cJ_h(U^*_h) = 0$, we see that
\bel{w1206e4}
\bal
\lt\| U^*_{h\t} - {\Pi_\t U^*_h}\rt\|_{L^2_\dbF(0,T; \dbL^2)}^2
&\leq \big[ \big(  D\widehat \cJ_h(U^*_h)-  D\widehat \cJ_{h}(\Pi_\t U^*_h),U^*_{h\t} - \Pi_\t U^*_h\big)_{L^2_\dbF(0,T; \dbL^2)} \\ 
& \q+ \big(  D\widehat \cJ_{h}(\Pi_\t U^*_h)- D\widehat \cJ_{h\t}(\Pi_\t U^*_h) ,U^*_{h\t} - \Pi_\t U^*_h
\big)_{L^2_\dbF(0,T; \dbL^2)}\big]\, .
\eal
\ee
Therefore,
\bel{w1212e1}
\bal
&\lt\| U^*_{h\t} - {\Pi_\t U^*_h}\rt\|_{L^2_\dbF(0,T; \dbL^2)}^2\\ 
 \leq& 3 \big[ \Vert  D\h\cJ_h(U^*_h)-  D\h\cJ_{h}(\Pi_\t U^*_h) \Vert_{L^2_\dbF(0,T;\dbL^2)}^2
+ \Vert \cT_h^1(\cS_h(\Pi_\t U^*_h))-Y_0(\cS_{h\t}(\Pi_\t U^*_h))
\Vert_{L^2_\dbF(0,T;\dbL^2)}^2\\ 
&\q+ \Vert Y_0(\cS_{h\t}(\Pi_\t U^*_h))-K_{h\t}\cS_{h\t}(\Pi_\t U^*_h)
\Vert_{L^2_\dbF(0,T;\dbL^2)}^2\big]\\ 
 =:& 3(I'+II'+III') \, .
\eal
\ee
We use \eqref{derivative-semidisc} and \eqref{pontr1a}
 to bound $I'$ as follows,
\bel{w1206e5}
\bal
I' =& \lt\| U^*_h - \Pi_\t U^*_h +   \cT_h^1 \big({\mathcal S}_h(\Pi_\t U^*_h) \big) 
-  \cT_h^1\big({\mathcal S}_h(U^*_h) \big)\rt\|_{L^2_\dbF(0,T; \dbL^2)}^2 \\
\leq&2 \big[\lt\| U^*_h - \Pi_\t U^*_h\rt\|_{L^2_\dbF(0,T; \dbL^2)}^2
  +\lt\|\cT_h^1 \big({\mathcal S}_h(\Pi_\t U^*_h) \big) -\cT_h^1\big({\mathcal S}_h(U^*_h) \big) \rt\|_{L^2_\dbF(0,T; \dbL^2)}^2\big]\, .
\eal
\ee
 By  stability properties of solutions to {\bf BSPDE}$_h$ \eqref{bshe1a} with $X^*_h=\cS_h(\Pi_\t U^*_h)$ and 
 {\bf SPDE$_h$} \eqref{w1013e1a} with $U_h=\Pi_\t U^*_h$, we obtain
\bel{w1206e6}
\bal
 & \lt\| \cT_h^1 \big({\mathcal S}_h(\Pi_\t U^*_h) \big) -\cT_h^1\big({\mathcal S}_h(U^*_h) \big) \rt\|_{L^2_\dbF(0,T; \dbL^2)}^2\\ 
 & \leq  C  \big[\lt\| \big(\cS_h(U^*_h)-\cS_h(\Pi_\t U^*_h)\big)(T) \rt\|^2_{L^2_{\mf_T}(\O;\dbL^2)}+ \lt\| \cS_h(U^*_h)-\cS_h(\Pi_\t U^*_h) \rt\|^2_{L^2_\dbF(0,T;\dbL^2)} \big] \\
 & \leq   C \me\big[\lt\| U^*_h-\Pi_\t U^*_h \rt\|^2_{L^2_\dbF(0,T;\dbL^2)}\big]\, .
\eal
\ee 
By the optimality condition \eqref{pontr1a}, estimates \eqref{w301e1}$_3$, \eqref{w301e1}$_5$ of Lemma \ref{w229l3},
we have
\bel{w1206e7}
\bal
 \lt\| U^*_h-\Pi_\t U^*_h \rt\|^2_{L^2_\dbF(0,T; \dbL^2)} 
 \leq C \lt\| Y_h-\Pi_\t Y_h \rt\|^2_{L^2_\dbF((0,T; \dbL^2)} 
  \leq C\t \|X_0\|_{\dbH_0^1}^2 \, . 
\eal
\ee

Next, we turn to $II'$. Triangular inequality leads to
\beq
II' 
\leq &2  \Big(\lt\| \cT_h^1 \big(\cS_{h}(\Pi_\t U^*_h) \big) -  \cT_h^1 \big(\cS_{h\t}(\Pi_\t U^*_h) \big) \rt\|^2_{L^2_\dbF(0,T; \dbL^2)} \\
  &+\lt\| \cT_h^1 \big(\cS_{h\t}(\Pi_\t U^*_h) \big) -Y_0 \big(\cS_{h\t}(\Pi_\t U^*_h)\big)  \rt\|^2_{L^2_\dbF(0,T; \dbL^2)}\Big)\\
=: &2 \big(II'_{1}+II'_{2}\big) \,. \\
\eeq
In order to bound $II'_{1}$, we use stability properties 
for {\bf SPDE}$_h$ \eqref{sde}, {\bf BSPDE}$_h$ \eqref{bshe1a}, in combination with the error estimate \eqref{euler1} for \eqref{sde} to conclude
$$II'_{1} \leq C \Big( \lt\Vert\lt( {\mathcal S}_h(\Pi_\t U^*_h) - {\mathcal S}_{h\t}(\Pi_\t U^*_h)\rt)(T)\rt\Vert_{L^2_{\mf_T}(\O; \dbL^2)}^2
+\Vert {\mathcal S}_h(\Pi_\t U^*_h) - {\mathcal S}_{h\t}(\Pi_\t U^*_h)\Vert_{L^2_\dbF(0,T; \dbL^2)}^2 \Big) 
 \leq C\t\, . $$
To bound $II'_{2}$, it is easy to see
\beq
\bal
II'_{2}\leq&  2\sum_{n=0}^{N-1}\me\Bigl[\int_{t_n}^{t_{n+1}}\lt\|Y_h(t;\cS_{h\t}(\Pi_\t U^*_h))-Y_h(t_n;\cS_{h\t}(\Pi_\t U^*_h))\rt\|^2_{\dbL^2} \rd t\Bigr] \\
 &+2T \max_{0\leq n\leq N}\me \bigl[\|Y_h(t_n;\cS_{h\t}(\Pi_\t U^*_h))-Y_0(t_n,\cS_{h\t}(\Pi_\t U^*_h))\|_{\dbL^2}^2\bigr]\, .
\eal
\eeq
Like \eqref{w301e2}, we can get
\bel{w310e1}
\bal
&\sum_{n=0}^{N-1}\me\Bigl[\int_{t_n}^{t_{n+1}}\lt\|Y_h(t;\cS_{h\t}(\Pi_\t U^*_h))-Y_h(t_n;\cS_{h\t}(\Pi_\t U^*_h))\rt\|^2_{\dbL^2} \rd t\Bigr] \\
& \leq C\t \me\Bigl[\int_0^T \|\D_h Y_h(t;\cS_{h\t}(\Pi_\t U^*_h))\|_{\dbL^2}^2
    +\|Z_h(t;\cS_{h\t}(\Pi_\t U^*_h))\|_{\dbL^2}^2
    +\|\cS_{h\t}(t;\Pi_\t U^*_h)\|_{\dbL^2}^2  \rd t\Bigr]\\
& \leq C\t\Big(\sup_{0\leq n\leq N}\me\bigl[\|\D_h\cS_{h\t}(t_n;\Pi_\t U^*_h)\|_{\dbL^2}^2\bigr]+\|\cS_{h\t}(\Pi_\t U^*_h)\|^2_{L^2_{\dbF}(0,T;\dbL^2)}\Big)\\
& \leq C\t   \lt[\|X_0\|_{\dbH_2}^2+\|\si\|^2_{C([0,T];\dbH^2)}\rt]\, .
\eal
\ee
Therefore,
by Lemma \ref{w308l1} and \eqref{w310e1},  we can get
\beq
II'_2\leq C\t \,.
\eeq

Finally,
Lemma \ref{w301l1} leads to 
\beq
III'\leq C\t .
\eeq

Now we insert above estimates into \eqref{w1212e1} to obtain assertion {\rm (i)}.

\ms

\no{\bf 2)}
For all $k=0,1,\cds,N$, we define $e_{X}^k=X^*_h(t_k)- X^*_{h\t}(t_k)$. Subtracting \eqref{w1003e7} from \eqref{w1013e1a} leads to
\beq
e_{X}^{k+1} - e_X^{k} 
=&  \t  \D_h e_X^{k+1} +  e_X^k\D_{k+1}W+
\t \big[U_h^*(t_k) - U_{h\tau}^*(t_k)\big] \\
&+ \int_{t_k}^{t_{k+1}}  \big( \big[X^*_h(s) - X^*_h(t_k)\big]+\big[\cR_h\si(s) - \cR_h\si(t_k)\big] \big)\, {\rm d}W(s)\\
&
+ \int_{t_k}^{t_{k+1}} \big( \D_h \big[X^*_h(s) - X^*_h(t_{k+1})\big] + \big[U^*_h(s) - U^*_h(t_k)\big]\big)
\, {\rm d}s\, .
\eeq
Testing with $e_{X}^{k+1}$, and using binomial formula, Poincar\'{e}'s inequality, independence, and absorption lead to
\beq
&\frac{1}{2} {\mathbb E}\big[ \Vert e_{X}^{k+1}\Vert^2_{\dbL^2} - \Vert e_{X}^{k}\Vert^2_{\dbL^2} + \frac{1}{2} \Vert e_{X}^{k+1} - e_{X}^{k}\Vert^2_{\dbL^2}\big] + \frac{\tau}{2} {\mathbb E} \bigl[\Vert \nabla e_{X}^{k+1} \Vert^2_{\dbL^2}\bigr]\\
 &\leq  \t \me\bigl[\|e^{k+1}_X \|_{\dbL^2}^2\bigr]+2 \t \me\bigl[\|e^{k}_X \|_{\dbL^2}^2\bigr]+ 
   \frac{\tau}{2} \me \bigl[\Vert U_h^*(t_k) - U_{h\tau}^*(t_k)\Vert^2_{\dbL^2}\bigr]  \\
&\q+ C \me \Bigl[\int_{t_k}^{t_{k+1}}  \Vert \nabla \big[ X^*_h(s) - X^*_h(t_{k+1})\big] \Vert^2_{\dbL^2}
+\Vert X^*_h(s) - X^*_h(t_k)\Vert_{\dbL^2}^2\\
&\qq\qq\qq+\Vert \si(s) - \si(t_k)\Vert_{\dbH_0^1}^2
+ \Vert U^*_h(s) - U^*_h(t_k)\Vert^2_{\dbL^2} \, {\rm d}s\Bigr]\, .
\eeq
By the discrete Gronwall's  inequality, and then taking the sum over all  $0\leq k\leq N-1$, and noting that $e_X^0=0$, 
we find that
\beq
&\max_{0\leq n\leq N}\me \bigl[\|e_X^n\|_{\dbL^2}^2\bigr] +\sum_{n=1}^N\t\me \bigl[\|\nb e_X^n\|_{\dbL^2}^2\bigr] \\
& \leq C\t \sum_{k=0}^{N-1}\me \bigl[\|U^*_h(t_k)-U^*_{h\t}(t_k)\|_{\dbL^2}^2\bigr] \\
&\q+C\sum_{k=0}^{N-1}\me \Bigl[\int_{t_k}^{t_{k+1}}   \Vert \nabla \big[ X^*_h(s) - X^*_h(t_{k+1})\big] \Vert^2_{\dbL^2}
+\Vert X^*_h(s) - X^*_h(t_k)\Vert_{\dbL^2}^2\\
&\qq\qq\qq+\Vert \si(s) - \si(t_k)\Vert_{\dbH_0^1}^2+ \Vert U^*_h(s) - U^*_h(t_k)\Vert^2_{\dbL^2}   \rd s \Bigr]\, .
\eeq
By \eqref{w1212e1}, the first term on the right-hand side is bounded by $C \t $. 
By \rf{w826e1}$_4$ and \rf{w1011e1}$_4$, we can bound the second and third terms by $C\t$.
The fourth term can be bounded by $C\t \|\si\|_{C^1([0,T];\dbH_0^1)}^2$.
By the optimal condition \eqref{pontr1a} and \eqref{w301e1}$_3$ in Lemma \ref{w229l3}, the last term is bounded  by $C\tau$.
That is assertion (ii).
\end{proof}

\section{The gradient descent method to solve problem {\bf SLQ}$_{h\t}$}\label{numopt}

By Theorem \ref{MP}, solving problem {\bf SLQ$_{h\t}$} is 
equivalent to solving the system of the coupled forward-backward difference equations 
\eqref{w1212e3} and \eqref{w1003e12}. We may exploit the variational
character of problem {\bf SLQ$_{h\t}$} to construct a gradient descent method
where approximate iterates of the optimal control $U^*_{h\t}$ in the 
Hilbert space $\dbU_{h\t}$ are obtained; see also \cite{Nesterov04,Kabanikhin12} for more details.
%
%

\ms

\begin{algorithm}\label{alg1}
Let $U_{h\t}^{(0)}\in \dbU_{h\t}$, and fix $\kappa > 0$. For any $\ell \in {\mathbb N}_0$, update $U_{h\t}^{(\ell)} \in {\mathbb U}_{h\tau}$ as follows:
\begin{enumerate}
\item[1.] Compute $X_{h\t}^{(\ell)}\in \dbX_{h\t}$ by 
\bel{w0115e4}
\left\{
\bal
& [\mathds{1} - \tau \Delta_h]X^{(\ell)}_{h\t}(t_{n+1})= X^{(\ell)}_{h\t}(t_n)+ \tau U^{(\ell)}_{h\t}(t_n) 
+ \big[X^{(\ell)}_{h\t}(t_n)+\cR_h\si(t_n)\big] \D_{n+1}W\\
&\qq\qq\qq\qq \q n=0,1,\cds,N-1\, ,\\
& X_{h\t}^{(\ell)}(0)=\cR_h X_0\, .
\eal
\right.
\ee
\item[2.] Use $X_{h\t}^{(\ell)}\in \dbX_{h\t}$ to compute $Y_{h\t}^{(\ell)}\in \dbX_{h\t}$ via
\begin{equation*}
\bal
Y_{h\t}^{(\ell)} (t_n)=&-\t \me\Big[ \sum_{j=n+1}^N A_0^{j-n} \prod_{k={n+2}}^j(1+ \D_kW) X^{(\ell)}_{h\t}(t_j)  \Big|\mf_{t_n} \Big]  \\
&-\a \me \Big[ A_0^{N-n} \prod_{k=n+2}^N(1+ \D_kW) X^{(\ell)}_{h\t}(T) \Big|\mf_{t_n}\Big] \q n=0,1,\cds,N-1\,.\\
\eal
\end{equation*}

\item[3.] Update $U_{h\t}^{(\ell+1)} \in {\mathbb U}_{h\tau}$ via
\beq
U^{(\ell+1)}_{h\t}=U^{(\ell)}_{h\t}-\frac 1 {\kappa} \big[U^{(\ell)}_{h\t} -Y_{h\t}^{(\ell)} \big] \, .
\eeq
\end{enumerate}
\end{algorithm}

If compared with \eqref{w1212e3}-\eqref{w1003e12}, steps 1 and 2 are now decoupled: the first step requires to solve a space-time discretization of  {\bf SPDE}$_h$ \rf{w212e3}$_1$, while the second requires to solve a space-time discretization of the {\bf BSPDE} \rf{w212e3}$_2$ which is not the numerical solution by the implicit Euler method (see Lemmas \ref{w229l1} and \ref{w301l1} for the difference).
A similar method to solve problem {\bf SLQ$_{h\t}$} 
has been proposed in \cite{Dunst-Prohl16, Prohl-Wang20}.
We refer to related works on how to approximate conditional expectations
 (e.g.~\cite{Bouchard-Touzi04, Gobet-Lemor-Warin05, Bender-Denk07, Wang-Zhang11, Wang-Wang-Lv-Zhang20}).


In the below,
we tend to present a lower bound for $\kappa$ and show convergence rate of Algorithm \ref{alg1}.
 For this purpose, we first recall Lipschitz continuity of 
$D\widehat{\mathcal J}_{h\tau}$.
%
%
By the definition of Fr\'echet derivative,
$$D^2\h\cJ_{h\t}(U_{h\t})= \big[\mathds{1} +L^* L+\a\h L^*\h L\big]U_{h\t}\,,$$
where operators $L,\,\h L$ are defined in \eqref{rep1}-\eqref{w1003e14}. Next, we can find $K:=\|\mathds{1} +L^* L+\a\h L^*\h L\|_{\cL(\dbU_{h\t};\dbU_{h\t})}$, such that
\beq
\|D\h\cJ_{h\t}\lt(U^1_{h\t}\rt)-D\h\cJ_{h\t}\lt(U^2_{h\t}\rt)\|_{\dbU_{h\t}}\leq K \lt\|U^1_{h\t}-U^2_{h\t}\rt\|_{\dbU_{h\t}}\, .
\eeq
Indeed,
noting that $\|(\mathds{1}-\t\D_h)^{-1}\|_{\cL(\dbU_{h\t},\dbU_{h\t})}\leq 1$,
we find that
\beq
\|LU_{h\t}\|_{\dbX_{h\t}}^2
=\sum_{n=1}^N \tau \me \| \t \sum_{j=0}^{n-1} \lt[(\mathds{1}-\t\D_h)^{-1}\rt]^{n-j}\prod_{k=j+2}^n\lt(1+\D_kW \rt) U_{h\t}(t_j) \|_{\dbL^2}^2  
\leq  T^2e^{T}\|U_{h\t}\|_{\dbU_{h\t}}^2\, .
\eeq
In a similar vein, we can prove that
\beq
\|\h LU_{h\t}\|_{L^2_{\mf_T}(\O;\dbL^2)}^2
\leq   Te^{T}\|U_{h\t}\|_{\dbU_{h\t}}^2\, .
\eeq
Hence
\beq
K=\|\mathds{1} +L^* L+\a\h L^*\h L\|_{\cL(\dbU_{h\t};\dbU_{h\t})}\leq 1+\a Te^{T}+T^2e^{T}\, .
\eeq
Since Algorithm \ref{alg1} is the gradient descent method for {\bf SLQ$_{h\t}$}, we have the following estimates.
\bt{gradient-rate} Suppose that $\kappa \geq K$.
Let $\ds \{U^{(\ell)}_{h\t}\}_{\ell \in {\mathbb N}_0} \subset {\mathbb U}_{h\tau}$ be  generated by Algorithm \ref{alg1}, and $U^*_{h\t}$ solve {\bf SLQ$_{h\t}$}. Then for $\ell=1,2,\cds$,
%
\begin{eqnarray*}
{\rm (i)} &&  \| U^{(\ell)}_{h\t}-U^*_{h\t} \|_{\dbU_{h\t}}^2\leq \Big(1-\frac 1 {\kappa}\Big)^{\ell}\| U^{(0)}_{h\t}-U^*_{h\t} \|_{\dbU_{h\t}}^2  \, ; \\ 
{\rm (ii)} &&\h\cJ_{h\t}(U^{(\ell)}_{h\t})-\h\cJ_{h\t}(U^*_{h\t})\leq \frac{2 \kappa \Vert U^{(0)}_{h\t}-U^*_{h\t} \Vert_{\dbU_{h\t}}^2}{\ell}\, ;\\
{\rm (iii)}&& \max_{0 \leq n \leq N}  {\mathbb E}\bigl[\|X_{h\t}^*(t_n)-X^{(\ell)}_{h\t}(t_n)\|_{\dbL^2}^2\bigr]
+\tau \sum_{n=1}^N {\mathbb E}   \bigl[\|X^*_{h\t}(t_n)- X^{(\ell)}_{h\t}(t_n)\|_{\dbH_0^1}^2\bigr] \\
&& \q
\leq C\Big(1-\frac 1 {\kappa}\Big)^{\ell}\| U^{(0)}_{h\t}-U^*_{h\t} \|_{\dbU_{h\t}}^2 \, ,
\end{eqnarray*}
where $C$ is independent of $h,\,\t,\,\ell$.
\et

\begin{proof}
{The estimates (i) and (ii) are standard for the gradient descent method (see e.g.~\cite[Theorem 1.2.4]{Nesterov04}); more details can also be found in \cite[Section 5]{Prohl-Wang20}. In the following, we restrict to assertion (iii).

For all $n=0,1,\cds,N$, define $\bar e_{X}^{n,\ell}=X^*_{h\t}(t_n)-X^{(\ell)}_{h\t}(t_n)$. Subtracting \rf{w0115e4} 
from \eqref{w1003e7} (where $U_{h\t}=U^*_{h\t}$) leads to
\beq
\bar e_{X}^{n+1,\ell} - \bar e_X^{n,\ell} 
=&  \t  \D_h \bar e_X^{n+1,\ell} +  \bar e_X^{n,\ell }\D_{n+1}W+
\t \big[U_{h\tau}^*(t_n)-U^{(\ell)}_{h\t}(t_n)\big] \, .
\eeq
Multiplication with $e_{X}^{n+1,\ell}$, and then taking expectations, as well as applying Cauchy-Schwartz inequality, we arrive at
\beq
&(1-\t) \me\big[ \Vert \bar e_{X}^{n+1,\ell}\Vert^2_{\dbL^2} \big] +2\t \me\big[\|\nb \bar e_X^{n+1,\ell}\|_{\dbL^2}^2 \big]  \\
 &\leq  (1+\t) \me\big[ \Vert \bar e_{X}^{n,\ell}\Vert^2_{\dbL^2} \big] +\t \me\big[\|U_{h\tau}^*(t_n)-U^{(\ell)}_{h\t}(t_n)\|_{\dbL^2}^2 \big]\, .
\eeq
By the discrete Gronwall's  inequality, then taking the sum over all  $0\leq n\leq N-1$, and utilizing the fact $\bar e_X^{0,\ell}=0$,
we find that
\beq
\max_{0\leq n\leq N}\me \bigl[\|\bar e_X^{n,\ell}\|_{\dbL^2}^2\bigr] +\sum_{n=1}^N\t\me \bigl[\|\nb \bar e_X^{n,\ell}\|_{\dbL^2}^2\bigr] 
 \leq C\t \sum_{n=0}^{N-1}\me \bigl[\|U_{h\tau}^*(t_n)-U^{(\ell)}_{h\t}(t_n)\|_{\dbL^2}^2\bigr] \, ,
\eeq
which, together with (i), imply assertion (iii).
}
\end{proof}

\section*{Acknowledgement}

This work was initiated when Yanqing Wang visited the University of T\"ubingen in 2019--2020, supported by a 
DAAD-K.C.~Wong Postdoctoral Fellowship.

%

\end{document}